\documentclass[a4paper,12pt]{amsart}
\usepackage[letterpaper, margin=1.2in]{geometry}
\usepackage{latexsym}
\usepackage{amssymb}
\usepackage{amsmath}
\usepackage{amsfonts}
 \usepackage{ulem}
\usepackage{tikz}
\usetikzlibrary{calc,arrows}
\usetikzlibrary{calc}
\usepackage{verbatim}
\newtheorem{thm}{Theorem}[section]

\newtheorem{lem}[thm]{Lemma}

\newtheorem{cor}[thm]{Corollary}

\theoremstyle{definition}

\newtheorem*{rem}{Remark}

\def\fph{\mathbb{F}_{\ph}}

\newcommand{\Z}{\mathbb Z}
\newcommand{\z}{\mathbb Z}
\newcommand{\Q}{\mathbb Q}

\newcommand{\F}{\mathbb F}

\newcommand{\q}{\mathbb Q}
\def\F{\mathbb{F}}

\newcommand{\p}{\mathfrak{p}}
\def\ol{\overline}
\def\al{\alpha}
\def\la{\lambda}

\def\om{\omega}
\def\th{\theta}
\def\md#1{\ \mbox{\rm(mod }{#1})}
\def\nph#1{N_{\ph}(#1)}
\def\npp#1{N_{\ph}^-(#1)}
\def\ph{\phi}
\def\ind{\text{ind}}

\newcounter{cs}
\stepcounter{cs}
\newcommand{\casos}{\begin{itemize}}
\newcommand{\fcasos}{\end{itemize}\setcounter{cs}{1}}

\newfont{\tit}{cmr12 scaled \magstep3}

\title[]{On  integral bases and monogenity of  pure octic number fields  with non-square free parameters}
\author{Lhoussain El Fadil and Istv\'an Ga\'al
}
\address{Faculty of Sciences Dhar El Mahraz, P.O. Box  1796 Atlas-Fes, Sidi Mohamed ben Abdellah University,  Morocco
}\email{lhouelfadil2@gmail.com}
\address{University of Debrecen, Mathematical Institute, H-4002 Debrecen Pf.400, Hungary    }
\email{gaal.istvan@unideb.hu }
\begin{document}
\date{}

\begin{abstract}
In all available papers, on power integral bases  of pure octic number fields $K$, generated by 
a root $\al$ of  a monic
irreducible polynomial $f(x)=x^8-m\in\Z[x]$, it was
assumed that $m\neq \pm 1$ is  square free. In this paper, we investigate the monogenity of any pure  octic number field, without the condition that $m$ is square free. 
We start by  calculating an integral basis of $\Z_K$, the ring  of integers of $K$. 
In particular, we characterize  when $\Z_K=\Z[\al]$.
We give  sufficient conditions on $m$, which guarantee that $K$ is not monogenic. We finish 
the paper by  investigating  the case when $m=a^u$, $u\in\{1,3,5,7\}$ 
and $a\neq \mp 1$ is a square free rational integer.
\end{abstract}
\maketitle
\section{Introduction}
Let $K=\Q(\al)$ be a  number field generated by a root  $\al$ of a monic irreducible polynomial 
$f(x)\in \Z[x]$ of degree $n$, denote by $\Z_K$ its ring of integers and $d_K$ its discriminant.
The ring $\Z_K$ is said to have a {\it power integral basis} if it has a $\Z$-basis $(1,\th,\cdots,\th^{n-1})$ 
for some $\th\in \Z_K$. That means $\Z_K=\Z[\th]$, that is  $\Z_K$ is mono-generated as a ring, with a single generator $\th$. In this case, the field $K$ is said to be {\it monogenic} and not monogenic otherwise.

The  problem of deciding the monogenity  of a number field  and constructing  
power integral bases is a classical problem of algebraic number theory, going back to
Dedekind \cite{R}, Hensel \cite{hensel08} and Hasse \cite{hasse63}.
This area is intensively studied even nowadays, cf. \cite{gaalbook} for the present state of this research.

\medskip

It is well-known that  $\Z_K$ is a free $\Z$-module of rank $n$. For any primitive
element $\th\in\Z_K$ (that is $K=\Q(\th)$), $\Z[\th]$ is a subgroup of $\Z_K$ of finite index.
We call
\[
\ind(\th)=(\Z_K:\Z[\th])
\]
the {\it index of $\th$}. 
Obviously, $\ind(\th)=1$ if and only if $(1,\th,\ldots,\th^{n-1})$ is an integral basis
of $\Z_K$. 
As it is known \cite{gaalbook}, we have
\[
\triangle(\th)=\ind(\th)^2\cdot d_K
\]
where $\triangle(\th)$ is the discriminant of $\th$.

The greatest common divisor $i(K)$ of the indices of the primitive elements $\th\in\Z_K$ of $K$
is called the {\it (field) index of $K$}.  
It is clear that if $i(K)>1$, then $K$ is not monogenic. On the other hand, $i(K)=1$ 
does not imply that $K$ is monogenic.
If a prime $p$ divides $i(K)$, 
then $p$ is called a prime {\it common index divisor} of $K$.
  
Let $(1,\omega_1,\ldots,\omega_{n-1})$ be an integral basis of $\Z_K$.
The discriminant $\triangle(L(X_1,\ldots,X_n))$ of the linear form 
$L(X_1,\ldots,X_{n-1})=\omega_1X_1+\ldots +\omega_{n-1}X_{n-1}$ can be written (cf. \cite{gaalbook}) as
\[
\triangle(L(X_1,\ldots,X_{n-1}))=(\ind(X_1,\ldots,X_{n-1}))^2\cdot d_K,
\]
where $\ind(X_1,\ldots,X_{n-1})$ is the {\it index form} corresponding to the 
integer basis $(1,\omega_1,\ldots,\omega_{n-1})$ having the property that for any
$\th=x_0+\omega_1x_1+\ldots +\omega_{n-1}x_{n-1}\in\Z_K$ (with $x_0,x_1,\ldots,x_{n-1}\in\Z$)
we have $\ind(\th)=|\ind(x_1,\ldots,x_{n-1})|$. Therefore $\th$ is 
a {\it generator of a power integral basis} if and only if
$x_1,\ldots,x_{n-1}\in\Z$ is a solution of the {\it index form equation}
$\ind(x_1,\ldots,x_{n-1})=\pm 1$.

According to our notation $K$ is generated by a root $\al$ of a monic irreducible polynomial $f\in\Z[x]$.
The index
\[
\ind (f)=(\Z_K:\Z[\alpha])
\]
is called the {\it index of the polynomial $f$} and we also have 
\[
\triangle(f)=\ind(f)^2\cdot d_K,
\]
$\triangle(f)$ denoting the discriminant of $f$.

\medskip

The main tool of our proofs will be the theory of {\it Newton polygons}.
If a prime $p$ does not divide $(\mathbb{Z}_K : \mathbb{Z}[\alpha])$, then a well-known theorem of Dedekind says that the  factorization of  $ p \mathbb{Z}_K$ can be derived directly  from the factorization of $\overline{f(x)}$ in $\F_p[x]$. In order to apply this  theorem in an effective way, one needs a
criterion to test  whether  $p$ divides   the index $(\Z_K:\Z[\al])$. 
Such a criterion was proved by Dedekind in 1878, (see  \cite[Theorem 6.1.4]{Co}). 
 When Dedekind's criterion fails, { then a method of Ore  $1928$,  based on {\it Newton polygon techniques}, can be used} in order to evaluate the  index $(\Z_K:\Z[\al])$, the
absolute discriminant of $K$, and the prime ideal factorization of the rational primes into powers of prime ideals of $\Z_K$ (see \cite{MN, O}).
{In case Ore's method also fails, then an algorithm developed by  Guardia, Montes, and Nart  \cite{GMN},
 based on   {\it higher order Newton polygons} can be used. Such an algorithm} gives after a finite number of iterations a complete answer on the index $(\Z_K:\Z[\al])$, the absolute discriminant $d_K$ of $K$, and the factorization of $p\Z_K$.

\medskip

Among all types of number fields, most investigations deal with pure number fields $K$ 
generated by a root of an irreducible polynomial $x^n-m$.  Assuming that $m$ is square free,
Ga\'al  and  Remete \cite{GR14} studied pure quartic fields,
Ahmad, Nakahara and Husnine \cite{NH, AN6}, Ahmad, Nakahara and Hameed \cite{AN}, El Fadil \cite{E6}
pure sextic fields. 
Applying the  index forms,  Ga\'al  and  Remete \cite{GR17} investigated
pure number fields of degrees $3\le n\le 9$.

The subject of our present paper is the monogenity of pure octic fields. 
For square free $m$ 
Hameed and Nakahara \cite{HN15}, proved that if $m \equiv 1 \; (\bmod \; 16)$, then the octic
number field generated by $m^{1/8}$ is not monogenic, but if $m \equiv 2, 3 \; (\bmod \; 4)$, then it is
monogenic.

In all above quoted  papers, the authors consider only pure octic number fields 
generated by a root of a monic polynomial $f(x)=x^8-m$, 
where $m\neq \pm 1$ is a square free rational integer.
Our purpose is to extend these results to arbitrary parameters $m$,
without assuming that $m$ is square free.
We start by calculating an  integral basis of $\Z_K$ in Theorem \ref{base}. We 
   give sufficient conditions on $m$, which guarantee the non monogenity   of $K$. We { conclude} the paper by studying the case where $m=a^u$, $u\in\{1,3,5,7\}$, and $a\neq \mp 1$ is a square free rational integer.

\medskip

\section{Main results}\label{main}

Throughout this paper, $m\neq \pm 1$ is a rational integer such that  the polynomial $f(x)=x^8-m$ is irreducible  over $\Q$. Let $\al$ be a root of $f(x)$, let $K=\Q(\al)$ with ring of integers $\Z_K$.
Replacing $\al$ by $\frac{\al}{p}$ and $m$ by $\frac{m}{p^8}$, and repeating this process until we get $\nu_p(m)<8$ for every  prime integer $p$, we can assume that $m =a_1a_2^2a_3^3a_4^4a_5^5a_6^6a_7^7$,
where $a_1,\dots,a_7$ are square free pairwise coprime
rational integers.
Set $A_2=a_4a_5a_6a_7$, $A_3=a_3a_4a_5a_6^2a_7^2$,  $A_4=a_2a_3a_4^2a_5^2a_6^3a_7^3$,   $A_5=a_2a_3a_4^2a_5^3a_6^3a_7^4$,  $A_6=a_2a_3^2a_4^3a_5^3a_6^4a_7^5$, and $A_7=a_2a_3^2a_4^3a_5^4a_6^5a_7^6$. 

\medskip 

The following theorem gives explicitly an  integral  basis of  $\z_K$. 

 \medskip 

\begin{thm}\label{base}
{
Using the above notations let 
$
m_2=\frac{m}{2^{\nu_2(m)}}
$ 
and let $u\in \Z$ such that $m_2u\equiv 1 \md{2^6}$.
In the following table $\bf{B}$ is an explicitly given integral  basis  of  $K$. }
 $$ \mbox{Table  A : } $$
	$$\begin{tabular}{| c| c| c| }
	\hline
	\mbox{conditions}& {\bf{B}}\\
	\hline
	$\nu_2(m)$\mbox{ is odd }& $(1,\al, \frac{\al^2}{A_2},\frac{\al^3}{A_3},\frac{\al^4}{A_4},\frac{\al^5}{A_5}, \frac{\al^6}{A_6}, \frac{\al^7}{A_7})$ \\
	\hline
$m\equiv 28\md{32}$&
$(1,\al, \frac{\al^2}{A_2},
\frac{\al^3}{A_3},\frac{\ph_2(\al)}{2A_4},\frac{\al\ph_2(\al)}{2A_5}, \frac{\al^2\ph_2(\al)}{4A_6}, \frac{\al^3\ph_2(\al)}{4A_7})$ \\
&{
$\ph_2(\al)=\al^4+2m_2u\al^2+2m_2u$} \\
		\hline
$m\equiv 12\md{32}$&
$(1,\al, \frac{\al^2}{A_2},
\frac{\al^3}{A_3},\frac{\ph_2(\al)}{2A_4},\frac{\al\ph_2(\al)}{2A_5}, \frac{\al^2\ph_2(\al)}{2A_6}, \frac{\al^3\ph_2(\al)}{4A_7})$ \\
&$\ph_2(\al)=\al^4+2m_2u\al^2+4m_2u\al+6m_2u$ \\
\hline
$m\equiv 4\md{16}$&
$(1,\al, \frac{\al^2}{A_2},
\frac{\al^3}{A_3},\frac{\ph_2(\al)}{2A_4},\frac{\al\ph_2(\al)}{2A_5}, \frac{\al^2\ph_2(\al)}{2A_6}, \frac{\al^3\ph_2(\al)}{2A_7})$ \\
&$\ph_2(\al)=\al^4+2m_2u$\\

\hline
$m\equiv 48\md{64}$&
{
$(1,\al, \frac{\al^2}{A_2},
\frac{\al^3}{A_3},\frac{\al^4}{A_4},\frac{\al\ph_2(\al)^2}{2A_5}, \frac{\al^6}{A_6}, \frac{\al^3\ph_2(\al)^2}{2A_7})$} \\
&$\ph_2(\al)=\al^2+2m_2u$ \\
\hline
$m\equiv 80\md{128}$&
$(1,\al, \frac{\al^2}{A_2},
\frac{\al\ph_2(\al)}{2A_3},\frac{(\ph_2(\al))^2}{2A_4},\frac{\al(\ph_2(\al))^2}{2A_5}, \frac{\al^2(\ph_2(\al))^2}{2A_6}, \frac{\al(\ph_2(\al))^3}{4A_7})$ \\
&$\ph_2(\al)=\al^2+2m_2u$\\
\hline
$m\equiv 144\md{256}$&
$(1,\al, \frac{\al^2}{A_2},
\frac{\al^3}{A_3},\frac{\ph_2(\al)^2}{2A_4},\frac{\al\ph_2(\al)^2}{2A_5}, \frac{\th}{4A_6}, \frac{\al\cdot \th}{4A_7})$ \\
&$\ph_2(\al)=\al^2+2m_2u$ \mbox{ and } $\th=\ph_2^3(\al)-8m_2u\ph_2^2(\al)+24m_2u\ph_2(\al)-32m_2u$\\
\hline
$m\equiv 16\md{256}$&
$(1,\al, \frac{\al^2}{A_2},
\frac{\al^3}{A_3},\frac{\ph_2(\al)^2}{2A_4},\frac{\al\ph_2(\al)^2}{2A_5}, \frac{\th}{4A_6}, \frac{\al\cdot \th}{8A_7})$ \\
&$\ph_2(\al)=\al^2+2m_2u$ \mbox{ and } $\th=\ph_2^3(\al)-8m_2u\ph_2^2(\al)+24m_2u\ph_2(\al)-32m_2u$\\

\hline
{
$m\equiv 448\md{512}$}&
$(1,\al, \frac{\al^2}{A_2},
\frac{\al^3}{A_3},\frac{\ph_2(\al)}{2A_4},\frac{\al\ph_2(\al)}{4A_5}, \frac{\al^2\ph_2(\al)}{2A_6}, \frac{\al^3\ph_2(\al)}{2A_7})$ \\
&$\ph_2(\al)=\al^4+4m_2u\al^3+12m_2u\al^2+24m_2u$ \\
\hline
{$m\equiv 192\md{512}$}&
$(1,\al, \frac{\al^2}{A_2},
\frac{\al^3}{A_3},\frac{\ph_2(\al)}{2A_4},\frac{\al\ph_2(\al)}{4A_5}, \frac{\al^2\ph_2(\al)}{4A_6}, \frac{\al^3\ph_2(\al)}{2A_7})$ \\
&{$\ph_2(\al)=\al^4+4m_2u\al^2+8m_2u$}  \\
\hline
$m\equiv 64\md{256}$&
$(1,\al, \frac{\al^2}{A_2},
\frac{\al^3}{A_3},\frac{\ph_2(\al)}{2A_4},\frac{\al\ph_2(\al)}{2A_5}, \frac{\al^2\ph_2(\al)}{2A_6}, \frac{\al^3\ph_2(\al)}{2A_7})$ \\
&{$\ph_2(\al)=\al^4+8m_2u$}\\
\hline
	$m\equiv 3\md4$ & $(1,\al, \frac{\al^2}{A_2},
\frac{\al^3}{A_3},\frac{\al^4}{A_4},\frac{\al^5}{A_5}, \frac{\al^6}{A_6}, \frac{\al^7}{A_7})$ \\
\hline
	$m\equiv 5\md8$& $(1,\al, \frac{\al^2}{A_2},
\frac{\al^3}{A_3},\frac{\al^4+m^4}{2A_4},\frac{\al^5+m^4\al}{2A_5}, \frac{\al^6+m^4\al^2}{2A_6}, \frac{\al^7+m^4\al^3}{2A_7})$ \\
\hline
	$m\equiv 9\md{16}$ &$(1,\al, \frac{\al^2}{A_2},
\frac{\al^3}{A_3},\frac{\al^4+m^4}{2A_4},\frac{\al^5+m^4\al}{2A_5}, 
\frac{\al^6-2m\al^5-m^2\al^4+m^2\al^2+2m\al+3m^2}{4A_6}$,
\\ 
&$\frac{\al^7-m\al^6+m^2\al^5-m\al^4+m^2\al^3-m\al^2+(m^2+4m)\al+m}{4A_7})$ \\
		\hline
$m\equiv 1\md{16}$ &$(1,\al, \frac{\al^2}{A_2},
\frac{\al^3}{A_3},\frac{\al^4+m^2}{2A_4},\frac{\al^5+m^2\al}{2A_5}, \frac{\al^6-2m\al^5-m^2\al^4+m^2\al^2+2m\al+3m^2}{4A_6}$, \\&$\frac{\al^7-m\al^6+m^2\al^5-m\al^4+m^2\al^3-m\al^2+(m^2+4m)\al+m}{8A_7})$ \\
\hline
	\end{tabular}$$
\end{thm}

\vspace{0.5cm}

As a consequence we obtain:

\begin{cor}\label{pib}
Let $K=\Q(\al)$ be a pure octic  number field generated by a  root $\al$ of a monic irreducible polynomial $f(x)=x^8-m\in \Z[x]$. Then { $\Z_K=\Z[\al]$ } if and only if $m$ is a square free integer and  $m\not\equiv 1 \md{4}$.
\end{cor}

\vspace{0.5cm}

In case of the integral basis 
$(1,\al, \frac{\al^2}{A_2},\frac{\al^3}{A_3},\frac{\al^4}{A_4},\frac{\al^5}{A_5}, \frac{\al^6}{A_6}, \frac{\al^7}{A_7})$,
considering the explicit form of factors of the index form we conclude:

\begin{thm} \label{index}
Keeping the notation of Section \ref{main}, if $\nu_2(m)$ is odd or $m\equiv 3\md4$, then
$$8a_1a_3a_5a_7|(a_2^2a_6^2\pm 1)$$
 is a necessary condition for the monogenity of $K$.
\end{thm}

By the $\pm$ sign we mean that for the monogenity of $K$ the divisibility condition must hold with 
at least one of the signs. Note that similar conditions can be derived also in the other cases of the integral basis,
but the calculation becomes far too complicated.

\vspace{0.5cm}

Our next main result  gives sufficient conditions on $m$ {for the non-monogenity  of $K$}. 
It relaxes the condition  $m$ is  square free required in  \cite{HN15, GR17}.
\begin{thm}\label{npib}
Let $K=\Q(\al)$ be a pure octic  number field generated by a root $\al$ of a monic irreducible polynomial $f(x)=x^8-m\in \Z[x]$. If one of the following conditions holds
\begin{enumerate}
\item
$m\equiv 1 \md{32}$, 
\item
 $m\equiv 272\md{512}$, 
\item
 $\nu_2(m)$ is odd and $a_2a_6\md8\in\{ 2,6 \}$,
\end{enumerate}
then $K$ is not monogenic.
\end{thm}

\vspace{0.5cm}

Finally, we consider monogenity of pure octic fields for $m$ of type $a^u$:

\begin{thm}\label{mono}
Assume that  $m=a^u$ with $a\neq \pm 1$   a square free rational integer and $u\in \{3,5,7\}$. Then 
\begin{enumerate}
\item 
If  $a\not\equiv
1 \md{4}$, then $K$ is monogenic
and $\Z_K$ is generated by  $\theta=\frac{\al^x}{a^y}$, where $(x,y)\in\Z^2$ is the unique  solution in non negative integers of the equation $ux-8y=1$ with $x<8$.
\item
If  $a\equiv
1 \md{4}$, then $K$ is not monogenic with the exception of  $a=-3$.
\end{enumerate}
\end{thm}

\vspace{0.5cm}

\section{{Preliminaries}}
In order to show  Theorem \ref{base} and  Theorem \ref{npib},  we recall some fundamental facts of Newton polygon techniques applied to algebraic number theory. Namely, the theorems on the index and on the prime ideal factorization. 
{For a detailed 
presentation of this theory we refer to the paper of Guardia and  Nart \cite{GN}.}

{Let $f(x)\in\Z[x]$ be a monic irreducible polynomial with a root $\al$,
let  $\ol{f(x)}=\prod_{i=1}^r \ol{\ph_i(x)}^{l_i}$ modulo $p$ be the factorization of $\ol{f(x)}$ into powers of monic irreducible coprime polynomials of $\F_p[x]$. Recall that  a well-known  theorem of  Dedekind says that: 
 \begin{thm}$($\cite[ Chapter I, Proposition 8.3]{Neu}$)$
 $$\mbox{If }   p \mbox{  does not divide the index } (\Z_K:\Z[\al]), \mbox{ then } p\Z_K=\prod_{i=1}^r \p_i^{l_i}, \mbox{ where } \p_i=p\Z_K+\phi_i(\al)\Z_K$$  and the residue degree of $\p_i$ is $f(\p_i)={\mbox{deg}}(\phi_i)$.
 \end{thm}
 In order to apply this  theorem, one needs a
criterion to test  whether  $p$ divides  the index $(\Z_K:\Z[\al])$.
  In $1878$, Dedekind  { proved } the following criterion:
 {
  \begin{thm}\label{Ded}$($Dedekind's criterion \cite[Theorem 6.1.4]{Co} and \cite{R}$)$\\
 For a number field $K$ generated by a  root $\al$ of a monic irreducible  polynomial $f(x)\in \Z[x]$ and a rational prime integer $p$, let $\overline{f}(x)=\prod_{i=1}^r\overline{\ph_i}^{l_i}(x)\md{p}$  be the factorization of   $\overline{f}(x)$ in $\F_p[x]$, where the polynomials $\ph_i\in\Z[x]$ are monic with their reductions irreducible over $\F_p$ and $\gcd(\overline{\ph_i},\overline{\ph_j})=1$ for every $i\neq j$. If we set
$M(x)=\cfrac{f(x)-\prod_{i=1}^r{\ph_i}^{l_i}(x)}{p}$, then $M(x)\in \Z[x]$ and the following statements are equivalent:
\begin{enumerate}
\item[1.]
$p$ does not divide the index $(\Z_K:\Z[\al])$.
\item[2.]
For every $i=1,\dots,r$, either $l_i=1$ or $l_i\ge 2$ and $\overline{\ph_i}(x)$ does not divide $\overline{M}(x)$ in $\F_p[x]$.
\end{enumerate}
\end{thm} }
 When Dedekind's criterion fails, that is,  $p$ divides the index $(\z_K:\z[\alpha])$  for every primitive element $\al\in \Z_K$ of $K$,{  then} it is not possible to obtain the prime ideal factorization of $p\z_K$ by{ Dedekind's theorem}.
In 1928, Ore developed   an alternative approach
for obtaining the index $(\z_K:\z[\alpha])$, the
absolute discriminant, and the prime ideal factorization of the rational primes in
a number field $K$ by using Newton polygons (see \cite{MN, O}). For more details on Newton polygon techniques, we refer to \cite{El, GMN}}.
 For any  prime integer  $p$,{ let $\nu_p$ be the $p$-adic valuation of $\Q$, $\Q_p$ its $p$-adic completion, and $\Z_p$ the ring of $p$-adic integers. Let  $\nu_p$  { be } the Gauss's extension of $\nu_p$ to $\Q_p(x)$; $\nu_p(P)=\mbox{min}(\nu_p(a_i), \, i=0,\dots,n)$ for any polynomial $P=\sum_{i=0}^na_ix^i\in\Q_p[x]$ and extended by $\nu_p(P/Q)=\nu_p(P)-\nu_p(Q)$ for every nonzero polynomials $P$ and $Q$ of $\Q_p[x]$.} Let 
$\phi\in\z_p[x]$ be a   monic polynomial  {\it whose reduction} is irreducible  in
$\F_p[x]$, let $\fph$ be 
the field $\frac{\F_p[x]}{(\overline{\phi})}$. For any
monic polynomial  $f(x)\in \z_p[x]$, upon  the Euclidean division
 by successive powers of $\ph$, we  expand $f(x)$ as follows:
$f(x)=\sum_{i=0}^la_i(x)\phi(x)^{i},$ called    the $\phi$-{\it expansion} of $f(x)$
 (for every $i$, deg$(a_i(x))<$
deg$(\phi)$). 
The $\ph$-{\it Newton polygon} of $f(x)$ with respect to $p$, is the lower boundary convex envelope of the set of points $\{(i,\nu_p(a_i(x))),\, a_i(x)\neq 0\}$ in the Euclidean plane, which we denote by $\nph{f}$.  The $\ph$-Newton polygon of $f$, is the process of joining the obtained edges  $S_1,\dots,S_r$ ordered by   increasing slopes, which  can be expressed as $\nph{f}=S_1+\dots + S_r$. 
For every side $S_i$ of $\nph{f}$, the length of $S_i$, denoted   $l(S_i)$  is the  length of its projection to the $x$-axis and its height , denoted   $h(S_i)$  is the  length of its projection to the $y$-axis. {Let $d(S_i)=\gcd(l(S_i), h(S_i))$ be the ramification degree of $S$.
 The {\it principal} $\ph$-{\it Newton polygon} of ${f}$},
 denoted $\npp{f}$, is the part of the  polygon $\nph{f}$, which is  determined by joining all sides of negative  slopes.
  For every side $S$ {of} $\npp{f}$, with initial point $(s, u_s)$ and length $l$, and for every 
$0\le i\le l$, we attach   the following
{{\ residue coefficient} $c_i\in\fph$ as follows:
$$c_{i}=
\displaystyle\left
\{\begin{array}{ll} 0,& \mbox{ if } (s+i,{\it u_{s+i}}) \mbox{ lies strictly
above } S,\\
\left(\dfrac{a_{s+i}(x)}{p^{{\it u_{s+i}}}}\right)
\,\,
\md{(p,\phi(x))},&\mbox{ if }(s+i,{\it u_{s+i}}) \mbox{ lies on }S,
\end{array}
\right.$$
where $(p,\phi(x))$ is the maximal ideal of $\z_p[x]$ generated by $p$ and $\ph$. 
Let $\la=-h/e$ be the slope of $S$, where  $h$ and $e$ are two positive coprime integers. Then  $d=l/e$ is the degree of $S$.  Notice that, 
the points  with integer coordinates lying on $S$ are exactly $$\displaystyle{(s,u_s),(s+e,u_{s}-h),\cdots, (s+de,u_{s}-dh)}.$$ Thus, if $i$ is not a multiple of $e$, then 
$(s+i, u_{s+i})$ does not lie in $S$, and so $c_i=0$. The polynomial
{$$f_S(y)=t_dy^d+t_{d-1}y^{d-1}+\cdots+t_{1}y+t_{0}\in\fph[y],$$}} { is } called  
the {\it residual polynomial} of $f(x)$ associated to the side $S$, where for every $i=0,\dots,d$,  $t_i=c_{ie}$.

    Let $\npp{f}=S_1+\dots + S_r$ be the principal $\ph$-Newton polygon of $f$ with respect to $p$.     We say that $f$ is a $\ph$-{\it regular polynomial} with respect to $p$, if  $f_{S_i}(y)$ is square free in $\fph[y]$ for every  $i=1,\dots,r$.      The polynomial $f$ is said to be  $p$-regular  if $\overline{f(x)}=\prod_{i=1}^r\overline{\ph_i}^{l_i}$ for some monic polynomials $\ph_1,\dots,\ph_t$ of $\Z[x]$ such that $\ol{\ph_1},\dots,\ol{\ph_t}$ are irreducible coprime polynomials over $\F_p$ and    $f$ is  a $\ph_i$-regular polynomial with respect to $p$ for every $i=1,\dots,t$.
 \smallskip
 
The  theorem of Ore plays  a fundamental role for proving our main Theorems:\\
  Let $\ph\in\Z_p[x]$ be a monic polynomial, with $\overline{\ph(x)}$ irreducible in $\F_p[x]$. As defined in \cite[Def. 1.3]{EMN},   the $\ph$-{\it index} of $f(x)$, denoted by $\ind_{\ph}(f)$, is  deg$(\ph)$ times the number of points with natural integer coordinates that lie below or on the polygon $\npp{f}$, strictly above the horizontal axis,{ and strictly beyond the vertical axis} (see $Figure\ 1$).
  
  \begin{figure}[htbp] 
 \centering
 \begin{tikzpicture}[x=1cm,y=0.5cm]
 \draw[latex-latex] (0,6) -- (0,0) -- (10,0) ;
 \draw[thick] (0,0) -- (-0.5,0);
 \draw[thick] (0,0) -- (0,-0.5);
 \node at (0,0) [below left,blue]{\footnotesize $0$};
 \draw[thick] plot coordinates{(0,5) (1,3) (5,1) (9,0)};
 \draw[thick, only marks, mark=x] plot coordinates{(1,1) (1,2) (1,3) (2,1)(2,2)     (3,1)  (3,2)  (4,1)(5,1)  };
 \node at (0.5,4.2) [above  ,blue]{\footnotesize $S_{1}$};
 \node at (3,2.2) [above   ,blue]{\footnotesize $S_{2}$};
 \node at (7,0.5) [above   ,blue]{\footnotesize $S_{3}$};
 \end{tikzpicture}
 \caption{    \large  $\npp{f}$.}
 \end{figure}
In the example of $Figure\ 1$, $\ind_\ph(f)=9\times$deg$(\ph)$.\\
\smallskip

  Now assume that $\overline{f(x)}=\prod_{i=1}^r\overline{\ph_i}^{l_i}$ is the factorization of $\overline{f(x)}$ in $\F_p[x]$, where every $\ph_i\in\Z[x]$ is monic polynomial, with $\overline{\ph_i(x)}$ is irreducible in $\F_p[x]$, $\overline{\ph_i(x)}$ and $\overline{\ph_j(x)}$ are coprime when $i\neq j$ and $i, j=1,\dots,t$.
For every $i=1,\dots,t$, let  $N_{\ph_i}^-(f)=S_{i1}+\dots+S_{ir_i}$ be the principal  $\ph_i$-Newton polygon of $f$ with respect to $p$. For every $j=1,\dots, r_i$,  let $f_{S_{ij}}(y)=\prod_{k=1}^{s_{ij}}\psi_{ijk}^{a_{ijk}}(y)$ be the factorization of $f_{S_{ij}}(y)$ in $\F_{\ph_i}[y]$. 
  Then we have the following index theorem of Ore (see \cite[Theorem 1.7 and Theorem 1.9]{EMN}, \cite[Theorem 3.9]{El}, {\cite[pp: 323--325]{MN}}, { and \cite{O}}).
	
\vspace{0.5cm}
	
 \begin{thm}\label{ore} $($Theorem of Ore$)$
 \begin{enumerate}
 \item
   $$\nu_p(\ind(f))\ge \sum_{i=1}^r \ind_{\ph_i}(f).$$  The equality holds if $f(x)$ is $p$-regular. 
\item
If  $f(x)$ is $p$-regular, then
$$p\Z_K=\prod_{i=1}^r\prod_{j=1}^{r_i}
\prod_{k=1}^{s_{ij}}\p^{e_{ij}}_{ijk},$$ 
 { is the factorization of $p\Z_K$ into powers of prime ideals of $\Z_K$ lying above $p$}, where {$e_{ij}=l_{ij}/d_{ij}$, $l_{ij}$ is the length of $S_{ij}$,  $d_{ij}$ is the ramification degree}
 of   $S_{ij}$, and $f_{ijk}=\mbox{deg}(\ph_i)\times \mbox{deg}(\psi_{ijk})$ is the residue degree of the prime ideal  $\p_{ijk}$ over $p$.
 \end{enumerate}
\end{thm}

\vspace{0.5cm}

{ When Ore's program fails;   that is if $f(x)$ is not $p$-regular, then it may happen that } some factors of $f(x)$ provided by Hensel’s  {lemma}  and refined by first order Newton polygon techniques are  not irreducible over $\Q_p$. {In this case} in order to complete the factorization of $f(x)$ in $\Q_p[x]$, Guardia, Montes, and Nart introduced the notion of  {\it high order Newton polygon}. They showed,  thanks to a theorem of index \cite[Theorem 4. 18]{GMN}, that  after a finite number of iterations this process yields all monic irreducible factors of  {$f(x)$ in $\Q_p[x]$}, all prime ideals of $\Z_K$ lying above a prime integer $p$, the  { $p$-valuation of the  index $(\Z_K:\Z[\al])$, and so up to a sign} the absolute discriminant of $K$. 

\vspace{0.5cm}

{ We recall here some fundamental techniques of Newton polygons of high order. For more details, we refer to \cite{GMN} and \cite{GN}. As introduced in \cite{GMN},  a {\it type of order} $r - 1$ } is  a data ${\bf{t}}  = (g_1(x), -\la_1, g_2(x), -\la_2,\dots, g_{r-1}(x), -\la_{r-1},\psi_{r-1}(x))$,
		where every $g_{i}(x)$ is a monic polynomial in $\z_p[x]$, $\la_i\in \q^+$, and
		$\psi_{r-1}(y)$ is a polynomial over a  finite field of {$\displaystyle p^{H}$ elements and $H=\displaystyle\prod_{i=0}^{r-2}f_i$},  with  $f_i=\mbox{deg}(\psi_{i}(x))$, satisfying the following recursive properties:
		\begin{enumerate}
			\item 
			$g_{1}(x)$ is irreducible modulo $p$,  $\psi_0(y) \in \F[y]$ ($\F_0=\F_p$) {is} the polynomial obtained
			by reduction of $g_{1}(x)$ modulo $p$, and  $\F_1 := \F_0[y]/(\psi_0(y)))$.
			\item
			For every $i=1,\dots,r-1$, the Newton polygon of $i^{th}$ order, $N_i(g_{i+1}(x))$,  has a single  sided of slope $-\la_i$.
			\item 
			For every $i=1,\dots,r-1$, the residual polynomial of $i^{th}$ order, $R_i(g_{i+1})(y)$ is an
			irreducible polynomial in $\F_i[y]$. Let $\psi_i(y)\in\F_i[y]$ be the monic polynomial
			determined by $R_i(g_{i+1})(y)\simeq \psi_i(y)$ (are equal up to multiplication by a nonzero element of $\F_i$), and  $\F_{i+1}= \F_i[y]/(\psi_i(y))$. Thus, $\F_0\subset \F_1\subset \dots\subset \F_r$ is a tower of finite fields.
			\item
			For every $i=1,\dots,r-1$, $g_{i+1}(x)$ has minimal degree among all monic polynomials
			in $\z_p[x]$ satisfying $(2)$ and $(3)$.
			\item 
			$\psi_{r-1}(y) \in\F_{r-1}[y]$ is a monic irreducible polynomial,  {$\psi_{r-1}(y)\neq y$}, and $\F_{r}= \F_{r-1}[Y]/(\psi_{r-1}(y))$.
		\end{enumerate}
		Here  the field
		$\F_i$ should not be confused with the finite field of $i$ elements.\\
		As for every $i=1,\dots,r-1$, the residual polynomial of the $i^{th}$ order, $R_i(g_{i+1})(y)$ is an
		irreducible polynomial in $\F_i[y]$, by  theorem of the product in order $i$, the polynomial $g_{i}(x)$ is irreducible in $\z_p[x]$. Let  $\om_0=[\nu_p,x,0]$  {be} the Gauss's extension of $\nu_p$ to $\q_p(x)$. As  for every $i=1,\dots,r-1$, the residual polynomial of the $i^{th}$ order, $R_i(g_{i+1})(y)$ is an		irreducible polynomial in $\F_i[y]$, then according to MacLane's notations  {and definitions} \cite{Mc},  $g_{i+1}(x)$ induces  a valuation  on $\q_p(x)$,  denoted by  $\om_{i+1}=e_i[\om_{i},g_{i+1},\la_{i+1}]$, where $\la_i=h_i/e_i$,   $e_i$ and  $h_i$ are  positive coprime integers. The valuation   $\om_{i+1}$ is called the {\it augmented valuation} of $\om_i$  with respect to $g_{i+1}$ and $\la_{i+1}$,  defined over   $\q_p[x]$  as follows:   
		$$
		\om_{i+1}(f(x))=\mbox{min}\{e_{i+1}\om_i(a_j^{i+1}(x))+jh_{i+1},\,j=0,\dots,n_{i+1}\},$$ 
		where $f(x)=\displaystyle\sum_{j=0}^{n_{i+1}}a_j^{i+1}(x)g_{i+1}^j(x)$ is the $g_{i+1}(x)$-expansion of $f(x)$. According to the terminology  in \cite{GMN}, the valuation  $\om_r$ is called the  $r^{th}$-order valuation associated to the data ${\bf{t}}$.
		For every order $r\ge 1$, the $g_r$-Newton polygon of $f(x)$, with respect  to  the valuation $\om_r$, { denoted $N_r(f)$ } is the lower boundary of the convex envelope of  the set of  points $\{(i,\mu_i), i=0,\dots, n_r\}$ in the  {Euclidean} plane, where $\mu_i=\om_r(a^r_i(x)g_r^i(x))$.  { Its principal part is denoted $N_r^-(f)$. }
		The following are the relevant theorems from Guardia-Montes-Nart's work (high order Newton polygon):

\vspace{0.5cm}

			\begin{thm}\label{highpol}$($\cite[Theorem 3.1]{GMN}$)$\\
				Let  {$f(x) \in \z_{p}[x]$ be a monic polynomial such that $\overline{f(x)}$ is a positive power of
				$\overline{\ph(x)}$ for some monic polynomial $\ph(x)$ such that $\overline{\ph(x)}$ is irreducible over $\F_0$}. If  $N_{r}^-(f)=S_1+\dots+S_g$ has $g$ sides, then we can
				split {$f(x)=F_0(x)\times F_1(x)\times\dots\times F_g(x)$} in $\z_p[x]$, such that $N_r(F_i)=S_i$ and
				$R_r({F_i})(y)=R_r({f})(y)$ up to  multiplication by a nonzero
				element of $\F_r$ for every $i=1,\cdots,g$.\\
							\end{thm}
			
\vspace{0.5cm}

			\begin{thm}\label{highres}$($\cite[Theorem 3.7]{GMN}$)$\\
				Let $f(x) \in \z_{p}[x]$ be a monic polynomial such that  $N_r^-(f)=S$ has a single side of finite slope $-\la_r$.
				If $R_r({f})(y)=\displaystyle\prod_{i=1}^t\psi_i(y)^{a_i}$ is the factorization in $\F_r[y]$, then  $f(x)$ splits as $f(x)=F_0(x)\times F_1(x)\times\cdots\times F_t(x)$ in $\z_p[x]$ such that $N_r({F_i})=S$ has  a single side of slope $-\la_r$ and $R_r({F_i})(y)=\psi_i(y)^{a_i}$  up to  multiplication by a nonzero
				element of $\F_r$ for every $i=1,\cdots,t$.
			\end{thm}
   \begin{rem}
{ The statement of  Theorem \ref{highpol} coincides with that given in \cite[Theorem 3.1]{GMN}, with $f_{\bf{t}}=F_1(x)\times\cdots\times F_g(x)$.\\
The statement of  Theorem \ref{highres} coincides with that given in   \cite[Theorem 3.7]{GMN}, with 
 $ f_{(\bf{t},\la_r)}=F_1(x)\times\cdots\times F_t(x)$. }
   \end{rem}

\vspace{0.5cm}

			{In \cite[Definition 4.11]{GMN}, the authors introduced a definition of the index of a polynomial $f(x)$ which is not necessarily irreducible over $\Q$ as follows: $\ind(f)=\sum_{i=1}^k\ind(F_i)+\displaystyle\sum_{1\le i<j\le k}\nu_p(res(F_i, F_j))$, where 	$f(x)=\prod_{i=1}^k F_i$ is the factorization of $f(x)$ in $\Q_p[x]$, $\ind(F_i)=\nu_p((\Z_i:\Z_p[\al_i]))$, $\al_i$ a root of $F_i(x)$ in ${\ol \Q_p}$, $\Z_i$ the integral closure of $\Z_p$ in $\Q_p(\al_i)$, and $res(F_i, F_j)$ is the resultant of $F_i$ and $F_j$. This definition of index extends the known one of $\ind(f)=(\Z_K:\Z[\al])$, with $\al\in\Z_K$  a primitive element of $K$ and $f$ its minimal polynomial over $\Q$.  For a fixed irreducible polynomial $F(x)$ in $\Z_p[x]$ and a fixed  data $${\bf{t}}  = (g_1(x), -\la_1, g_2(x), -\la_2,\dots, g_{r-1}(x), -\la_{r-1},\psi_{r-1}(x)),$$   the authors introduced
   in \cite[Definition 4.15]{GMN}, the notion of $r^{th}$-{\it order index of  $F(x)\in \Z_p[x]$} as follows: Let $N_r(F)$ be the Newton polygon of $r^{th}$-order with respect to the data  {\bf{t}} and   $\ind_r(F)=f_0\cdots f_{r-1}\ind(N_r(f))$, where   $f_i=\mbox{deg}(\psi_{i}(x))$ and $\ind(N_r(F))$ is the index of the polygon $N_r(F)$; the number of points with natural integer coordinates that lie below or on the polygon $N_r(F)$, strictly above the horizontal line of equation $y=\om_r(F)$, { and strictly beyond the vertical axis}}. 
		They showed the following theorem on the index  which generalizes the theorem of index of Ore as follows:
  {  Let $F(x)\in \Z[x]$ be a monic polynomial, { which is irreducible over $\Q_p$} and $${\bf{t}}  = (g_1(x), -\la_1, g_2(x), -\la_2,\dots, g_{r-1}(x), -\la_{r-1},\psi_{r-1}(x))$$ a fixed data such that $N_r(F)$ has a single side of negative slope ({ $\psi_{r-1}(y)$ divides $R_r(F)(y)$}). Then we have the following theorem of index:}
  
			\begin{thm}\label{thmind}$($\cite[Theorem 4.18]{GMN}$)$\\
			$$ \nu_p(\ind(F))\ge  \ind_1(F) +\dots + \ind_r(F).$$ 
The equality holds if and only if $\ind_{r+1}(F) = 0$. 
			\end{thm}
			Recall that by definition $\ind(N_{r+1}(F))=0$ if and only if $N_{r+1}(F)$ has a single side of length $1$ or height $1$. By  \cite[Lemma 2.17]{GMN} $(2)$, if $R_r(F)$ is square free, then the length of $N_r(F)$
 is $1$. Thus if $R_r(F)$ is square free, then $\ind_{r+1}(F)=0$, and so the equality { $ \nu_p(\ind(F))=  \ind_1(F) +\dots + \ind_r(F)$} holds. \\

  {
In order to complete the calculation of the index of any separable polynomial $f\in \Z[x]$, the  authors  introduced an iterative method to evaluate $\nu_p(res(F_i,F_j)$, with  $f(x)=\prod_{i=1}^k F_i$  the factorization of $f(x)$ in $\Q_p[x]$, as follows: if ${\ol {F_i}}$ and ${\ol {F_j}}$ are coprime modulo $p$, then  $res_r(F_i,F_j)=0$. If  ${\ol {F_i}}$ and ${\ol {F_j}}$ are congruent to a power of a monic irreducible polynomial ${\ol {g_1}}\in \F_0[x]$, then  the  authors  introduced the $r^{th}$-resultant of $F_i$ and $F_i$ as follows:
$$  res_r(F_i,F_j) = f_0\cdots f_{r-1}\times {\min}( E_iH_j, E_jH_i),$$
where $E_i$ and $H_i$ are the length and height of the sides $S_i$ of 
$N^{-}_r(F_i)$. Recall the following convention: If $S_i$ is reduced to a single point, then  $E_i=H_i=0$ and $R_r(F_i)$ is a constant of $\F_r$.\\ 
 Thanks to \cite[Lemma 4.8]{GMN} and \cite[Theorem 4.10]{GMN}, we have   the following: $$ \nu_p(res(F_i,F_j))\ge  res_1(F_i,F_j) +\dots + res_r(F_i,F_j).$$  
{ Moreover, if the data  $${\bf{t}}  = (g_1(x), -\la_1, g_2(x), -\la_2,\dots, g_{r-1}(x), -\la_{r-1},\psi_{r-1}(x))$$ satisfies the condition $N_r(F_i)$ has a single side of negative slope;  $\psi_r(y)$ divides $R_r(F_i)(y)$, then  the equality holds if and only if $R_{r+1}(F_i)(y)$ and $R_{r+1}(F_j)(y)$ are coprime}.\\
{ In particular, if ${\ol {F_i}}$ and ${\ol {F_j}}$ are coprime modulo $p$, then  $\nu_p(res(F_i,F_j))=0$. \\
  If for some integer $r=1,\dots, k$, $\psi_r(y)$ does not divide $R_r(F_s)$ for some $s=1,2$, then the equality holds.}
\vspace{0.5cm}

\section{Proofs of main results}

\vspace{0.5cm}

\noindent
{\it Proof of Theorem \ref{base}}.\\
 {During this proof, $\F_i$ is the $i^{th}$ field of  the tower provided by Montes algorithm.} Since  $\triangle(f)=\mp 8^8m^7$ is the discriminant of $f$, thanks to the  formula {$\triangle(f)=(\Z_K:\Z[\al])^2d_K$}, the prime candidates to divide the index $(\Z_K:\Z[\al])$ are those dividing $2\cdot m$.
 \begin{enumerate}
\item
 Let $p$ be a prime dividing $m$. Then $\ol{f(x)}=x^8\md{p}$. For $\ph=x$, $\nph{f}=S$ has a single side joining the points $(0,\nu_p(m))$  
 and $(8,0)$. Then {
$\ind_1(f)=\ind_\ph(f)=\nu_p(\prod_{i=2}^7A_i)$.}  Also by \cite[Theorem 2.7]{EMN}, {
$(1,\al, \frac{\al^2}{A_2},
\frac{\al^3}{A_3},\frac{\al^4}{A_4},\frac{\al^5}{A_5} ,\frac{\al^6}{A_6}, \frac{\al^7}{A_7})$} is a free $\Z$-sub-module of  $\Z_K$.
Now, let $d=${gcd}$(\nu_p(m),8)$. Then $R_1(f)(y)=y^d-\ol{m_p}$, where $m_p=m/p^{\nu_p(m)}$. If $2$ does not divide $\nu_p(m)$, then $d=1$, $R_1(f)(y)$ is irreducible over $\F_1$, and so by Theorem \ref{ore}, {$\ind(f)=\ind_\ph(f)$} and $(1,\al, \frac{\al^2}{A_2},
\frac{\al^3}{A_3},\frac{\al^4}{A_4},\frac{\al^5}{A_5}, \frac{\al^6}{A_6}, \frac{\al^7}{A_7})$ is a $p$-integral basis of $\Z_K$.  If $2$ divides $d$, then $d\in\{2,4\}$. In this case if  $p\neq 2$,  then $R_1(f)(y)=y^d-\ol{m_p}$ is  square free over $\F_1\simeq \F_0$ (because $\ph=x$, and so  {deg$(\ph)=1=[\F_1:\F_0]$}). Thus by Theorem \ref{ore}, {$\nu_p(\ind(f))=\ind_1(f)$} and 
$(1,\al, \frac{\al^2}{A_2},
\frac{\al^3}{A_3},\frac{\al^4}{A_4},\frac{\al^5}{A_5}, \frac{\al^6}{A_6}, \frac{\al^7}{A_7})$ is a $p$-integral basis of $\Z_K$. It follows that if $2$ divides $m$ and $\nu_2(m)$ is odd, then  $(1,\al, \frac{\al^2}{A_2},
\frac{\al^3}{A_3},\frac{\al^4}{A_4},\frac{\al^5}{A_5} ,\frac{\al^6}{A_6} ,\frac{\al^7}{A_7})$ is an integral basis of $\Z_K$.
\\
If $d\in\{2,4\}$ and $p=2$, then $R_1(f)(y)=(y-\ol{m_2})^d$ is not square free and we have to use second order Newton polygon techniques.
 \begin{enumerate}
 \item
If $\nu_2(m)=2$, then for  $\ph=x$, we have $\ol{f(x)}=\ph^8\md2$, $N_1(f)=S$ has a single side of slope $-\la_1=-1/4$, $e_1=4$, and $R_1(f)(y)=y^2+1=(y+1)^2$. According to the definitions and notations of \cite{GMN, Mc},  let $\om_2=e_1[\nu_2,\la_1]$ be the  valuation of second order Newton polygon defined by $\om_2(a)=e_1\nu_2(a)=4\nu_2(a)$ for every $a\in \Q_2$ and $\om_2(x)=e_1\la_1=1$. Let $\ph_2=x^4+2$ be {a key polynomial of $\om_2$} and 
$f(x)=\ph_2^2-4\ph_2+(4-m)$  the $\ph_2$-expansion of $f(x)$. As $\om_2(\ph_2)=4$, $\om_2(\ph_2^2)=8$, and $\om_2(4\ph_2)=12$. It follows that:
\begin{enumerate}
 \item
If $\nu_2(4-m)=3$; {$m\equiv 12 \md{16}$}, then there are $2$ cases :\\
 If {$\nu_2(12-m)= 4$; $m\equiv 28\md{32}$}, then for $\ph_2=x^4+2x^2+2$, we have $f(x)=\ph_2^2-4x^2\ph_2-(4+m)$ is  the $\ph_2$-expansion of $f(x)$. As $\om_2(\ph_2^2)=8$,  $\om_2((-4x^2)\ph_2)=14$, and  $\om_2(4+m)\ge 20$, we conclude that {  if $\om_2(4+m)= 20$, then  $N_2(f)$ has a single side,  $\ind_2(f)=6$, and
  $R_2(f)(y)=y^2+y+1$. { Thus $f(x)$ is irreducible over $\Q_2$ and $\nu_2(\ind(f))=\ind_1(f)+\ind_2(f)=4+6=10$}. If $\om_2(4+m)> 20$, then  $N_2(f)$ has two sides of degree $1$ each.
   { $f(x)=F_1(x)F_2(x)$ and $\nu_2(\ind(f))=\nu_2(\ind(F_1))+\nu_2(\ind(F_2))+\nu_2(res(F_1,F_2))=0+0+res_1(F_1,F_2)+res_2(F_1,F_2)=4+6=10$}.  Based on the  polygon $N_1(f)$, we conclude that $V(\al)=1/4$. Similarly, based on  $N_2(f)$, we conclude that
$V(\ph_2(\al))\ge 5/2$, {and so $(1,\al, {\al^2},
{\al^3},\frac{\ph_2(\al)}{4},\frac{\al \ph_2(\al)}{4}, \frac{\al^2\ph_2(\al)}{8},\frac{\al^3\ph_2(\al)}{8})$ is a $2$-integral basis of $\Z_K$. Since the Montes algorithm is local; the algorithm  {provides $p$-integral bases},  sometimes we have to replace $\ph_2(x)$ by an $\omega_2$-equivalent polynomial. For example, in our case, if $\nu_p(A_i)\ge 1$ for some $i=2,\dots,7$ and for some odd prime integer $p$, then $\frac{\al^3\ph_2(\al)}{8A_7}$ is not $p$-integral ($V(\frac{\al^3\ph_2(\al)}{8A_7})<0$ for some valuation $V$ of $K$ extending $\nu_p$). So, we have to replace $\ph_2(x)$ by $g(x)=x^4+2m_2ux^2+2m_2u$ with  $u$ an integer which satisfies $um_2\equiv 1 \md{32}$} and  show that  $\mathcal{B}=(1,\al, \frac{\al^2}{A_2},
\frac{\al^3}{A_3},\frac{g(\al)}{2A_4},\frac{\al g(\al)}{2A_5}, \frac{\al^2g(\al)}{4A_6},\frac{\al^3g(\al)}{4A_7})$ is an integral basis of $\Z_K$. Since for every prime integer $p$, $\nu_p(\Z_K:\Z[\al])=\nu_p(2^6\prod_{i=2}^7A_i)$, we need only to show that every element of $\bf{B}$ is integral over $\Z$. By the definition of $g(x)$ and by the first point of this proof the $V$-valuation of each element of $\bf{B}$ is greater or equal than $0$ for every valuation $V$ of $K$ extending $\nu_p$ for every odd prime integer $p$. Let us show the same result for $p=2$. For this reason we need to   give a lower bound of $V(\al)$ and $V(\ph_2(\al))$ for every valuation $V$ of $K$ extending $\nu_2$. Let $V$ be a valuation of $K$ extending $\nu_2$. Since $V(\al)=1/4$, $V(\ph_2(\al))\ge 5/2$, $\nu_2(g(x)-\ph_2(x))\ge 5$, and $\al\in \Z_K$, we conclude that   $V(g(\al)-\ph_2(\al))\ge 5$. Thus $V(g(\al))\ge 5/2$, and so by a simple verification, the $V$-valuation of each element of $\bf{B}$ is greater or equal than $0$.  Hence      $(1,\al, \frac{\al^2}{A_2},
\frac{\al^3}{A_3},\frac{\ph_2(\al)}{2A_4},\frac{\al\ph_2(\al)}{2A_5}, \frac{\al^2\ph_2(\al)}{4A_6},\frac{\al^3\ph_2(\al)}{4A_7})$ is an integral basis of $\Z_K$.\\
In the remainder of this proof, these techniques will be repeated. So in every case, we give an adequate $\ph_2(x)$ for which $f(x)$ is regular with respect to $\om_2$, we give the $\ph_2$-expansion of $f(x)$, and a lower bound of $V(\ph_2(\al))$ for every valuation $V$ of $K$ extending $\nu_2$.}  
\\
 If {$\nu_2(12-m)\ge 5$ ($m\equiv 12\md{32}$)}, then for $\ph_2=x^4+2x^2+4x+6$,
  $f(x)=\ph_2^2+(-8-4x^2-8x)\ph_2+(12-m+16x^3+32x^2+32x)$ is  the $\ph_2$-expansion of $f(x)$.
 As $\om_2(\ph_2^2)=8$,  $\om_2((-8-4x^2-8x)\ph_2)=14$, and 
 $\om_2(12-m+16x^3+32x^2+32x)=\om_2(16x^3)=19$, we conclude that $N_2(f)=S$ has a single side joining $(0,19)$ and $(2,8)$. Therefore $\ind_2(f)=5$, the side $S$ is of degree $1$, $f(x)$ is irreducible over $\Q_2$, and   by Theorem \ref{thmind} $\nu_2(\ind(f))=\ind_1(f)+5$. Thus $(1,\al, \frac{\al^2}{A_2},
\frac{\al^3}{A_3},\frac{\ph_2(\al)}{2A_4},\frac{\al\ph_2(\al)}{2A_5}, \frac{\al^2\ph_2(\al)}{2A_6},\frac{\al^3\ph_2(\al)}{4A_7})$ is a $2$-integral basis of $\Z_K$.
  Now by replacing $\ph_2(x)$ by $g(x)=x^4+2m_2ux^2+4m_2ux+6m_2u$ with  $u$ an integer which satisfies $um_2\equiv 1 \md{16}$, we conclude that  $\mathcal{B}=(1,\al, \frac{\al^2}{A_2},
\frac{\al^3}{A_3},\frac{g(\al)}{2A_4},\frac{\al g(\al)}{2A_5}, \frac{\al^2g(\al)}{2A_6},\frac{\al^3g(\al)}{4A_7})$ is an integral basis of $\Z_K$.  
\item
If $\nu_2(4-m)= 4$;  $m\equiv 20 \md{32}$, then  $\om_2(4-m)=16$ and  $N_2(f)=T$ has a single side joining the points $(0,16)$, $(1,12)$, and $(2,8)$.  Thus, $\ind_1(f)=4$ and $R_2(f)(y)=y^2+y+1$ is irreducible over $\F_2=\F_1=\F_0$ (because deg$(\ph)=$deg$(\psi_1)=1$). Hence , $f(x)$ is irreducible over $\Q_2$ and by Theorem \ref{thmind}
 $\nu_2(\ind(f))=\ind_1(f)+4$. Based on the polygon $N_2(f)$, we conclude that $V(\ph_2(\al))=2$. Replacing $\ph_2(x)$ by $x^4+2m_2u$ with an integer $u$ satisfying $um_2\equiv 1 \md{16}$, we conclude that       $(1,\al, \frac{\al^2}{A_2},
\frac{\al^3}{A_3},\frac{\ph_2(\al)}{2A_4},\frac{\al\ph_2(\al)}{2A_5}, \frac{\al^2\ph_2(\al)}{2A_6},\frac{\al^3\ph_2(\al)}{2A_7})$ is an integral basis of $\Z_K$.
\item
If $\nu_2(4-m)\ge 5$; $m\equiv 4\md{32}$, then  $\om_2(4-m)\ge 20$ and  $N_2(f)=T_1+T_2$ has two sides joining $(0,v)$, $(1,12)$, and $(2,8)$ with $v\ge 20$. Thus each side is of  degree $1$, and so $f(x)=F_1(x)F_2(x)$ with every $F_i(x)$  irreducible over $\Q_2$. Thus
 $\nu_2(\ind(f))=\ind_1(F_1)+\ind_1(F_2)+res_1(F_1,F_2)+\ind_2(F_1)+\ind_2(F_2)+res_2(F_1,F_2)=res_1(F_1,F_2)+res_2(F_1,F_2)=4+4=8$. 
 {Based on $N_2(f)$, we have $V(\ph_2(\al))\ge 2$, and so by replacing $\ph_2(x)$ by $x^4+2m_2u$ with  an integer $u$ satisfying $um_2\equiv 1 \md{16}$,
   we conclude that \\ $(1,\al, \frac{\al^2}{A_2},
\frac{\al^3}{A_3},\frac{\al^4+2m_2u}{2A_4},\frac{\al^5+2m_2u \al}{2A_5} ,\frac{\al^6+2m_2u\al^2}{2A_6} ,\frac{\al^7+2m_2u\al^3}{2A_7})$} is an integral basis of $\Z_K$. \\
 {Note that the two cases $(i)$ and $(ii)$ could 
combined into one case, namely $\nu_2(4-m)\ge 4$.}
\end{enumerate}
\item
If $\nu_2(m)=4$; $m\equiv 16 \md{32}$, then for $\ph=x$,  $\ol{f(x)}=\ph^8\md2$, $\nph{f}=S$ has a single side of slope $-1/2$, $R_1(f)(y)=y^4+1=(y+1)^4$. Let $\om_2$ be the valuation of second order Newton polygon defined by $\om_2(a)=2\nu_2(a)$ for every $a\in \Q_2$ and $\om_2(x)=1$. Let $\ph_2=x^2+2$ and $f(x)=\ph_2^4-8\ph_2^3+24\ph_2^2-32\ph_2+(16-m)$  the $\ph_2$-expansion of $f(x)$.  {Since $m\equiv 16 \md{32}$,} then $\nu_2(m-16)\ge 5$. It follows that:
\begin{enumerate}
 \item
 If $\nu_2(16-m)=5$; $m\equiv 48\md{64}$, then  $N_2(f)=T$ has a single side joining the points $(0,10)$ and $(4,8)$ with slope $-1/2$, $\ind_2(f)=2$, and residual polynomial $R_2(f)(y)=(y+1)^2$. Let us use the third order Newton polygon associated to the data $t=(x, 1/2,\ph_2, 1/2, \ph_3)$, where $\ph_3=\ph_2^2(x)+4x=x^4+4x^2+4x+4$ is a key polynomial of $\om_2$ and let  $\om_3=2[\om_2,1/2]$ be the valuation of third order Newton polygon; $\om_3(a)=4\nu_2(a)$
 for every $a\in \Q_2$, $\om_3(x)=2$,  and $\om_3(\ph_2)=2(2+1/2)=5$. Let 
 $f(x)=\ph_3^2+(8 - 8 x^2-8x )\ph_3+(-m -48-16x^2-32x+32x^3)$ be the {$\ph_3$}-expansion of $f(x)$.   Since $m =48+64k$ for some integer $k$, $-m -48-16x^2-32x+32x^3=-16\ph_2(x)-64+64k-32x+32x^3$, and so $\om_3(-m -48-16x^2-32x+32x^3)=\om_3(16)+\om_3(\ph_2)=21$. Thus $N_3(f)=T$ has a single side joining the points $(0,21)$ and $(2,20)$. It follows that its height is $1$,  $\ind_3(f)=0$ and $f(x)$ is irreducible over $\Q_2$. By Theorem \ref{thmind} $\nu_2(\ind(f))=\ind_1(f)+\ind_2(f)=\ind_1(f)+2$. Based on $N_1(f)$, we have  $V(\al)=1/2$ and based on $N_2(f)$, we get   $V(\ph_2(\al))=5/4$. Replacing $\ph_2(x)$ by by $x^2+2m_2u$ with $um_2\equiv 1 \md{16}$, we get  $(1,\al, \frac{\al^2}{A_2},
\frac{\al^3}{A_3},\frac{\al^4}{A_4},\frac{\al\ph_2^2(\al)}{2A_5}, \frac{\al^6}{A_6},\frac{\al^3\ph_2^2(\al)}{2A_7})$ is an integral basis of $\Z_K$.
\item
 If  {$\nu_2(16-m)=6$ ($m\equiv 80\md{128}$)}, then  $N_2(f)=T_2$ has a single side joining the points $(0,12)$ and $(4,8)$. Thus $\ind_2(f)=6$. As $T$ is of degree $4$ and { $R_2(f)(y)=y^4+y^2+1=(y^2+y+1)^2$}, we have to use third order Newton polygon techniques. {Let $\om_3$ be the valuation of third order Newton polygon, $\ph_3=\ph_2^2+2x\ph_2+4x^2$ and $f(x)=\ph_3^2+(-(4x+12)\ph_2+16+24x)\ph_3+(-(128x+32)\ph_2+192x-m+80)$. Since $\om_3(\ph_2)=3$, we conclude that $N_3(f)=T_3$ has a single side joining the points $(0,12)$ and $(2,13)$, which  is of degree $1$. Thus  $\ind_3(f)=0$, $f(x)$ is irreducible over $\Q_2$ and by Theorem \ref{thmind}   $\nu_2(\ind(f))=\ind_1(f)+6$.
 Based on $N_1(f)$ and $N_2(f)$, we get $V(\al)=1/2$ and   $V(\ph_2(\al))=3/2$. 
  Replacing $\ph_2(x)$ by $x^2+2m_2u$ by $um_2\equiv 1 \md{16}$, $V(\ph_2(\al))=3/2$ and  $(1,\al, \frac{\al^2}{A_2},
\frac{\al\ph_2(\al)}{2A_3},\frac{\ph_2^2(\al)}{2A_4},\frac{\al\ph_2^2(\al)}{2A_5}, \frac{\al^2\ph_2^2(\al)}{2A_6},\frac{\al\ph_2^3(\al)}{4A_7})$} is an integral basis of $\Z_K$.
\item
  If $\nu_2(16-m)=  7$; $m\equiv 144\md{256}$, then $\om_2(16-m)=  14$ and $N_2(f)=T_1+T_2$   has $2$ sides joining the points $(0,14)$, $(1,12)$, $(2,10)$, and $(4,8)$. Thus  { $\ind_2(f)=7$.  Since the  attached residual polynomials of $f(x)$ are $R_{12}(f)(y)=y^2+y+1$ and $R_{22}(f)(y)=(y+1)^2$,} in order to complete the calculation of the index $\nu_2(\ind(f))$, we have to use third order Newton polygon. Let $\ph_3=\ph_2+4x=x^2+4x+2$ be the key polynomial of $\om_2$,  $\om_3=[\om_2,1]$ the valuation of third order Newton polygon; $\om_3(a)=2\nu_2(a)$
 for every $a\in \Q_2$, $\om_3(x)=1$, and $\om_3(\ph_2)=3$. Let 
 $f(x)=\ph_3^4+(88-16x)\ph_3^3+(728-544x)\ph_3^2+(96-3776x)\ph_3-({6528x}+m+3824)$ be the $\ph_3$-expansion of $f(x)$.
 As {$\om_3((728-544x)\ph_3^2)=12$, $\om_3((96-3776x)\ph_3)=13$}, and $\om_3(6528x+m+3824)=13$, $N_3^-(f)=T$ has a single side joining the points $(0,13)$ and $(2,12)$. Thus $T$ is of height $1$, $\ind_3(f)=0$ and $f(x)$ is irreducible over $\Q_2$. By Theorem \ref{thmind}  { $\nu_2(\ind(f))=\ind_1(f)+\ind_2(f)=\ind_1(f)+7$.
 Based on $N_1(f)$, $N_2(f)$, and $N_3(f)$, we get $V(\al)=1/2$ and $V(\ph_2(\al))\ge 3/2$. 
Let $\th=\ph_2^3(\al)-8\ph_2^2(\al)+24\ph_2(\al)-32$. We need to show that $\frac{\th}{2^5}\in \Z_K$. For this reason we have to show that $V(\th)\ge 5$ for every valuation $V$ of $K$ extending $\nu_2$.  Since $N_2(f)$ has two sides of slopes $-2$ and $-1$, by Theorem \ref{highpol}, $f(x)=f_1(x)\times f_2(x)$ in $\Z_2[x]$  and there are two distinct valuations $V_1$ and $V_2$ of $K$ extending $\nu_2$ which satisfy $V_1(\ph_2(\al))=2$ and $V_2(\ph_2(\al))=3/2$.  If $V(\ph_2(\al))=2$, then a simple verification shows that $V(\th)\ge 0$. If $V(\ph_2(\al))=3/2$, as $f(\al)=0$, we have $\th=\frac{m-16}{\ph_2(\al)}$, and so $V(\th)=7-3/2\ge 5$. Hence $\frac{\th}{2^5}\in\Z_K$ and  
 $(1,\al, \frac{\al^2)}{2},
\frac{\al\ph_2(\al)}{4},\frac{\ph_2^2(\al)}{8},\frac{\al\ph_2^2(\al)}{8}, \frac{\th}{32},\frac{\al\th}{32})$ is a $2$-integral basis of $\Z_K$.
By replacing $\ph_2(x)$ by $x^2+2m_2u$ with $m_2u\equiv 1\md{32}$, we get   $(1,\al, \frac{\al^2}{A_2},
\frac{\al\ph_2(\al)}{2A_3},\frac{\ph_2^2(\al)}{2A_4},\frac{\al\ph_2^2(\al)}{2A_5}, \frac{\th}{4A_6},\frac{\al\th}{4A_7})$ is an integral basis of $\Z_K$, where $\th=\ph_2^3(\al)-8m_2u\ph_2^2(\al)+24m_2u\ph_2(\al)-32m_2u$.}\\
\item
If $\nu_2(16-m)\ge 8$; {$m\equiv 16\md{256}$}, then let $v=\om_2(16-m)$, { $\ph_2=x^2-2x+2$ and $f(x)=\ph_2^4+(16+8x)\ph_2^3+(-40+32x)\ph_2^2-32x\ph_2+(16-m)$ be the $\ph_2$-expansion of $f(x)$. If $v\ge 9$, then $N_2(f)=T_1+T_2+T_3$  has $3$ sides joining the points $(0,2v)$, $(1,13)$, $(2,10)$, and {$(4,8)$}. In this case $f(x)=F_1(x)\times F_2(x)\times F_3(x)$ with  $F_1(x)$ and $F_2(x)$ irreducible over $\Q_2$. 
For $F_3(x)$, since   its  residual polynomial attached to $T_3$ is $R_{2}(f)(y)=(y+1)^2$, we have to use third order Newton polygon techniques. Let $\ph_3=\ph_2+2x=x^2+4x+2$ be the key polynomial of $\om_2$,  $\om_3$ the valuation of third order Newton polygon and
 $f(x)=\ph_3^4+(88-16x)\ph_3^3+(728-544x)\ph_3^2+(96-3776x)\ph_3-(6528x+m+3824)$  the $\ph_3$-expansion of $f(x)$.
 As {$\om_3((728-544x)\ph_3^2)=12$, $\om_3((96-3776x)\ph_3)=13$} and $\om_3(6528x+m+3824)=13$, $N_3^-(f)=T$ has a single side joining the points $(2,13)$ and $(2,12)$ with height $1$. Thus $F_3(x)$ is irreducible over $\Q_2$ and $\ind_3(f)=0$, $\nu_2(\ind(f))=\ind_1(f)+\ind_2(f)=\ind_1(f)+8$.
Thus  $\nu_2(\ind(f))=\ind_1(F_1)+\ind_1(F_2)+\ind_1(F_3)+ 
res_1(F_1,F_2)+res_1(F_1,F_3)+res_1(F_2,F_3)
+\ind_2(F_1)+\ind_2(F_2)+\ind_2(F_3)
+ 
res_2(F_1,F_2)+res_2(F_1,F_3)+res_2(F_2,F_3)
=0+0+1+1+2+2+0+0+2+2+2+2=\ind_1(f)+8$.
 If $v=8$, then $N_2(f)=T_1+T_2$  has $2$ sides joining the points $(0,16)$, $(1,13)$, $(2,10)$, and {$(4,8)$}. 
 In this case $f(x)=F_1(x)\times F_2(x)$. Since   the residual  polynomial of $f(x)$ attached to $T_1$; $R_{2}(f)(y)=y^2+y+1$ is irreducible over $\F_2$, we conclude that  $F_1(x)$ is irreducible over $\Q_2$. 
For $F_2(x)$, Since   the attached residual polynomial of the last side is $R_{2}(f)(y)=(y+1)^2$, we have to use third order Newton polygon techniques. Let $\ph_3=\ph_2+2x=x^2+4x+2$ be the key polynomial of $\om_2$,  $\om_3$ the valuation of third order Newton polygon, and
 $f(x)=\ph_3^4+(88-16x)\ph_3^3+(728-544x)\ph_3^2+(96-3776x)\ph_3-(6528x+m+3824)$  the $\ph_3$-expansion of $f(x)$.
 As {$\om_3((728-544x)\ph_3^2)=12$, $\om_3((96-3776x)\ph_3)=13$}, and $\om_3(6528x+m+3824)=13$, $N_3^-(f)=T$ has a single side joining the points $(2,13)$ and $(2,12)$ with height $1$. In this case we have also $\nu_2(\ind(f))=\ind_1(f)+8$.  Let $V$ be a valuation of $K$ extending $\nu_2$. In both cases, based on $N_1(f)$ and $N_2(f)$, we get $V(\al)=1/2$ and  $V(\ph_2(\al))\ge 3/2$. 
 If $V(\ph_2(\al))\ge 5/2$, then $V(\ph_2^3(\al)+(16+8\al)\ph_2^2(\al)+(-40+32\al)\ph_2(\al)-32\al)\ge 11/2$, and so
 $V(\frac{\al\th}{2^6}\ge 0$, where $\th=\ph_2^3(\al)+(16+8\al)\ph_2^2(\al)+(-40+32\al)\ph_2(\al)-32\al$.
 If $V(\ph_2(\al))=3/2$, then $V(\frac{\ph_2^3(\al)+(16+8\al)\ph_2^2(\al)+(-40+32\al)\ph_2(\al)-32\al}{2^5})=V(\frac{16-m}{2^5\ph_2(\al)})\ge 8-5$. Thus $V(\frac{\th}{2^5})\ge 0$ and $V(\frac{\al\th}{2^6})\ge 0$ for every valuation $V$ of $K$ extending $\nu_2$. Thus 
 $(1,\al, \frac{\al^2}{2},
\frac{\al\ph_2(\al)}{4},\frac{\ph_2^2(\al)}{8},\frac{\al\ph_2^2(\al)}{8}, \frac{\th}{32},\frac{\al\th}{64})$ is a $2$-integral basis of $\Z_K$.
 Replace $\ph_2(x)$ by $x^2-2m_2ux+2m_2u$ with $m_2u\equiv 1\md{64}$, we conclude that 
 $(1,\al, \frac{\al^2}{A_2},
\frac{\al\ph_2(\al)}{2A_3},\frac{\ph_2^2(\al)}{2A_4},\frac{\al\ph_2^2(\al)}{2A_5}, \frac{\th}{4A_6},\frac{\al\th}{8A_7})$ is an integral basis of $\Z_K$, where $\th=\ph_2^3(\al)+(16m_2u+8m_2u\al)\ph_2^2(\al)+(-40m_2u+32m_2u\al)\ph_2(\al)-32m_2u\al$.}
 \end{enumerate}
  \item
If $\nu_2(m)=6$; {$m\equiv 64\md{128}$}, then for $\ph=x$, we have $\ol{f(x)}=\ph^8\md2$, $\nph{f}=S$ has a single side of slope $-\la_1=-3/4$, $e_1=4$, and $R_1(f)(y)=y^2+1=(y+1)^2$. Let $\om_2$ be the  valuation of second order Newton polygon defined by $\om_2(a)=e_1\nu_2(a)=4\nu_2(a)$ for every $a\in \Q_2$ and $\om_2(x)=e_1\la_1=3$. It follows that:
\begin{enumerate}
\item
If  { $\nu_2(m-64)=7$ ($m\equiv 192\md{256}$), then there are two case:}\\
If { $\nu_2(m-192)=8$ ($m\equiv 448\md{512}$)}, then for $\ph_2=x^4+4x^3+12x^2+24$,
 we have $f(x)=\ph_2^2+
 (-224-8x^3-8x^2+32x)\ph_2+(4800-m+512x^3+2304x^2-768x)$ is  the $\ph_2$-expansion of $f(x)$. As  $\om_2(\ph_2^2)=24$,   $\om_2((-224-8x^3-8x^2+32x)\ph_2)=30$, and $\om_2(4800-m+512x^3+2304x^2-768x)=\om_2(768x)=35$, we conclude that $N_2(f)$ has   a single  side of degree $1$. Thus $f(x)$ is irreducible over $\Q_2$,  $\ind_2(f)=5$ and $\nu_2(\ind(f))=\ind_1(f)+5$. Since $V(\al)=3/4$ and $V(\ph_2(\al))= 35/8$. By replacing $\ph_2(x)$ by $x^4+4um_2x^3+12um_2x^2+24um_2$ with $m_2u\equiv 1\md{32}$, we get   $(1,\al, \frac{\al^2}{A_2},
\frac{\al^3}{A_3},\frac{\ph_2(\al)}{2A_4},\frac{\al\ph_2(\al)}{4A_5}, \frac{\al^2\ph_2(\al)}{2A_6},\frac{\al^3\ph_2(\al)}{2A_7})$ is an integral basis of $\Z_K$.\\
If  { $\nu_2(m-192)\ge 9$ ($m\equiv 192\md{512}$)}, then for $\ph_2=x^4+4x^2+8$,  we have $f(x)=\ph_2^2-8x^2\ph_2-(64+m)$ is  the $\ph_2$-expansion of $f(x)$. As  $\om_2(\ph_2^2)=24$,   $\om_2(-8x^2\ph_2)=30$, and $\om_2(64+m)=\om_2(m-192)\ge 36$, if  $\nu_2(m-192)=9$ ($\om_2(m-192)= 36$), then $N_2(f)$ has a single side with $R_2(f)(y)=y^2+y+1$, 
$f(x)$ is irreducible over $\Q_2$ and  $\nu_2(\ind(f))=\ind_1(f)+6$.

If  $\nu_2(m-192)\ge 10$, then $N_2(f)$ has 
$2$ sides of degree $1$ each. Thus $f(x)=F_1(x)\times F_2(x)$ with each $F_i$ is irreducible over $\Q_2$. In this case we have also  $\nu_2(\ind(f))=\ind_1(f)+6$. Since $V(\ph_2(\al))\ge 9/2$ and  $V(\al)=3/4$, if we replace $\ph_2(x)$ by $x^4+4m_2ux+8m_2u$ with $um_2\equiv 1\md{32}$, then $(1,\al, \frac{\al^2}{A_2},
\frac{\al^3}{A_3},\frac{\ph_2(\al)}{2A_4},\frac{\al\ph_2(\al)}{4A_5}, \frac{\al^2\ph_2(\al)}{4A_6},\frac{\al^3\ph_2(\al)}{2A_7})$ is an integral basis of $\Z_K$.
\item
If {$\nu_2(m-64)\ge 8$ ( $m\equiv 64\md{256}$)}, then for $\ph_2=x^4+8$,  we have $f(x)=\ph_2^2-16\ph_2+(64-m)$ is  the $\ph_2$-expansion of $f(x)$. Since  $\om_2(\ph_2^2)=24$,   $\om_2(16\ph_2)=28$, and $\om_2(64-m)\ge 32$, it follows that if $\om_2(64-m)= 32$, then  $N_2(f)$ has a single side with $R_2(f)(y)=y^2+y+1$. So $\nu_2(\ind(f))=\ind_2(f)+4$. If $\om_2(64-m)> 32$, then
$N_2(f)$ has two  side with degree $1$ each. In this case, we have also  $\nu_2(\ind(f))=\ind_2(f)+4$. In both cases we have $V(\al)=3/2$ and $V(\ph_2(\al))\ge 4$. Replacing $\ph_2$ by  $\ph_2=x^4+8um_2$ with $um_2\equiv 1 \md{32}$, we get   $(1,\al, \frac{\al^2}{A_2},
\frac{\al^3}{A_3},\frac{\ph_2(\al)}{2A_4},\frac{\al\ph_2(\al)}{2A_5}, \frac{\al^2\ph_2(\al)}{2A_6},\frac{\al^3\ph_2(\al)}{2A_7})$  an integral basis of $\Z_K$.
\end{enumerate}
\end{enumerate}
\item
 If  $m\equiv 1 \md{2}$, then $f(x)\equiv
(x-1)^8\md{2}$. Let $\ph=x-1$,  $f(x)=\ph^8+8\ph^7+28\ph^6+56\ph^5+70\ph^4+56\ph^3+28\ph^2+8\ph+1-m$. 
 It follows that:
 \begin{enumerate}
 \item
If { $\nu_2(1-m)=1$}; $m\equiv 3\md4$, then $\nph{f}=S$ has a single side joining $(0,1)$ and $(8,0)$. Thus $\nu_2(\ind(f))=\ind_{\ph}(f)=0$. It follows
 that  $(1,\al, \frac{\al^2}{A_2},
\frac{\al\ph_2(\al)}{A_3},\frac{\al^4}{A_4},\frac{\al^5}{A_5}, \frac{\al^6}{A_6},\frac{\al^7}{A_7})$ is an integral basis of $\Z_K$.   
\item
If   { $\nu_2(1-m)=2$};  $m\equiv 5 \md{8}$, then $\nph{f}=S_{1}$ has a single  side  joining $(0,2)$, $(4,1)$ and $(8,0)$ with residual polynomial $R_1(f)(y)=y^2+y+1$, which is irreducible over $\F_1=\F_0$. Thus $\nu_2(\ind(f))=\ind_{1}(f)=4$ and \\
 $(1,\al, \frac{\al^2}{A_2},
\frac{\al^3}{A_3},\frac{\al^4+m^4}{2A_4},\frac{\al^5+m^4\al}{2A_5}, \frac{\al^6+m^4\al^2}{2A_6}, \frac{\al^7+m^4\al^3}{2A_7})$ is an integral basis of $\Z_K$.
\item
If  { $\nu_2(1-m)=3$};   $m\equiv 9 \md{16}$, then $\nph{f}=S_{1}+S_{2}$ has $2$ sides  joining $(0,3)$, $(2,2)$, $(4,1)$, and $(8,0)$. Thus
$R_{11}(f)(y)=y^2+y+1$ and  $R_{21}(f)(y)=y+1$. It follows by Theorem \ref{ore} that $\nu_2(\ind(f))=\ind_{1}(f)=6$ and \\ $(1,\al, \frac{\al^2}{A_2},
\frac{\al^3}{A_3},\frac{\al^4+m^4}{2A_4},\frac{\al^5+m^4\al}{2A_5}, \frac{\al^6+2m\al^5+3m^2\al^4+m\al^2+2m\al+3m^2}{4A_6},$ \\ $ \frac{\al^7+2m\al^6+3m^2\al^5+m\al^3+2m\al^2+3m^2\al}{4A_7})$  is an integral basis of $\Z_K$.
\item
If    { $\nu_2(1-m)=4$};  $m\equiv 17 \md{32}$, then $\nph{f}=S_{1}+S_{2}+S_3$ has $3$ sides  joining $(0,4)$, $(1,3)$, $(2,2)$, $(4,1)$, and $(8,0)$. Thus $R_{11}(f)(y)=y^2+y+1$ and  $R_{i1}(f)(y)=y+1$  for every $i=2,3$. It follows by Theorem \ref{ore} that $\nu_2(\ind(f))=\ind_{1}(f)=7$  
and \\ $(1,\al, \frac{\al^2}{A_2},
\frac{\al^3}{A_3},\frac{\al^4+m^4}{2A_4},\frac{\al^5+m^4\al}{2A_5}, \frac{\al^6-2m\al^5-m^2\al^4+m^2\al^2+2m\al+3m^2}{4A_6},\\ \frac{\al^7-m\al^6+m^2\al^5-m\al^4+m^2\al^3-m\al^2+(m^2+4m)\al+m}{8A_7})$  is an integral basis of $\Z_K$.
\item
If    { $\nu_2(1-m)\ge 5$}; $m\equiv 1 \md{32}$, then $\nph{f}=S_{1}+S_{2}+S_3+S_4$ has $4$ sides  joining $(0,v)$, $(1,3)$, $(2,2)$, $(4,1)$, and $(8,0)$ with $v=\nu_2(1-m)\ge 5$. Thus $R_{i1}(f)(y)=y+1$  for every $i=1,2,3,4$. It follows by Theorem \ref{ore} that $\nu_2(\ind(f))=\ind_{1}(f)=7$ 
and \\ $(1,\al, \frac{\al^2}{A_2},
\frac{\al^3}{A_3},\frac{\al^4+m^4}{2A_4},\frac{\al^5+m^4\al}{2A_5}, \frac{\al^6-2m\al^5-m^2\al^4+m^2\al^2+2m\al+3m^2}{4A_6},\\
\frac{\al^7-m\al^6+m^2\al^5-m\al^4+m^2\al^3-m\al^2+(m^2+4m)\al+m}{8A_7})$  is an integral basis of $\Z_K$.
\end{enumerate}
\end{enumerate}
\hfill$\Box$

\vspace{0.5cm}

\noindent
{\it Proof of Theorem \ref{index}.}\\
Let $\alpha$ be a root of $f(x)=x^8-m$ as above, $\zeta$ a primitive eighth root of unity.
{The conjugates of $\alpha$ are
$\alpha^{(1)}=\alpha,\alpha^{(2)}=\zeta\alpha,\ldots,\alpha^{(7)}=\zeta^6\alpha,\alpha^{(8)}=\zeta^7\alpha$.}
{{For every $k=1,\dots,8$, let $\sigma_k$ be the embedding of $K$ defined by  $\sigma_k(\al)=\zeta^k \al$. Then  $\sigma_1, \sigma_2,\dots,\sigma_8$ are the embeddings of the number field  $K$}.}

Let $L(x)=\alpha x_1+\alpha^2 x_2+\alpha^3 x_3+\alpha^4 x_4+\alpha^5 x_5+\alpha^6 x_6+\alpha^7 x_7$
and denote by $L^{(i)}(x)$ its conjugates, corresponding to $\alpha^{(i)} \;\; (1\le i\le 8)$. 
Let 
\[
H_{ij}(x)=\frac{L^{(i)}(x)-L^{(j)}(x)}{\alpha^{(i)}-\alpha^{(j)}}
\]
for $1\le i<j\le 8$. 
{
We have $\sigma(H_{ij}(x))=H_{\overline{i+1},\overline{j+1}}(x)$
where $\overline{i+1}=i+1$ for $1\leq i\le 7$ and $\overline{i+1}=1$ for $i=8$, and similarly  $\overline{j+1}$.
}
 Set
\begin{eqnarray*}
F_4(x)&=&H_{15}(x)\cdot H_{26}(x)\cdot H_{37}(x)\cdot H_{48}(x),\\
F_8(x)&=&H_{12}(x)\cdot H_{23}(x)\cdot H_{34}(x)\cdot H_{45}(x)
\cdot H_{56}(x)\cdot H_{67}(x)\cdot H_{78}(x)\cdot H_{81}(x).
\end{eqnarray*}
{
Using symmetric polynomials these factors were explicitly calculated
by using Maple \cite{maple} which showed that they have rational integer coefficients.
(It also follows from the fact that these products remain fixed under any element of the Galois group of the normal
closure of $K$.)
}
These polynomials divide the index form 
\[
J(x_1,\ldots,x_7)=\prod_{1\le i<j\le 8}|H_{ij}(x)|
\]
of the basis $(1,\alpha,\alpha^2,\alpha^3,\alpha^4,\alpha^5,\alpha^6,\alpha^7)$ of $K$,
which is also of integer coefficients. Denote by
$F_{16}(x)$ the third factor of the above product. ($F_{16}(x)$ is of degree $16$, since the index form has 28 factors.) 

Therefore we obtain
\[
J(x_1,\ldots,x_7)=
\pm F_4(x)\cdot F_8(x)\cdot F_{16}(x),
\]
where the factors $F_4(x), F_8(x), F_{16}(x)$ have integer coefficients, 
of degrees $4$, $8$, $16$, 
having $10$, $169$, $9185$ terms, respectively.

If $\nu_2(m)$ is odd or $m\equiv 3\md4$, then $K=\Q(\sqrt[8]{m})$
has integral basis 
\[
\left(1,\alpha,\frac{\alpha^2}{A_2},\frac{\alpha^3}{A_3},\frac{\alpha^4}{A_4},
\frac{\alpha^5}{A_5},\frac{\alpha^6}{A_6},\frac{\alpha^7}{A_7}\right).
\]
Therefore 
\[
\sqrt{|D_K|}=\frac{\sqrt{|D(\alpha)|}}{A_2A_3A_4A_5A_6A_7}
=\frac{\prod_{1\le i<j\le 8} |\alpha^{(i)}-\alpha^{(j)}|}{A_2A_3A_4A_5A_6A_7}.
\]
Any algebraic integer $\vartheta\in \Z_K$ can be written as
\[
\vartheta=\sum_{k=0}^7 y_k\frac{\alpha^k}{A_k}
=\frac{1}{A_7}\sum_{k=0}^7\frac{A_7}{A_k}y_k\alpha^k=\frac{1}{A_7}\sum_{k=0}^7 z_k\alpha^k,
\]
where $A_0=1$, $A_7/A_k$ and $y_k$ are integers, $z_k=\frac{A_7}{A_k}y_k,\;\; (k=0,1,\ldots,7)$.
Hence
\[
\ind(\vartheta)=\frac{1}{\sqrt{|D_K|}}\prod_{1\le i<j\le 8}|\vartheta^{(i)}-\vartheta^{(j)}|
\]
\[
=
\frac{A_2A_3A_4A_5A_6A_7}{\prod_{1\le i<j\le 8} |\alpha^{(i)}-\alpha^{(j)}|}\cdot \frac{1}{A_7^{28}}\cdot 
\prod_{1\le i<j\le 8}(L^{(i)}(z)-L^{(j)}(z))
\]
\[
=
\frac{A_2A_3A_4A_5A_6A_7}{A_7^{28}}
\prod_{1\le i<j\le 8}\frac{L^{(i)}(z)-L^{(j)}(z)}{\alpha^{(i)}-\alpha^{(j)}}.
\]
We obtain
\[
\ind(\vartheta)=\frac{A_2A_3A_4A_5A_6A_7}{A_7^{28}}\cdot F_4(z)\cdot F_8(z)\cdot F_{16}(z).
\]
In $F_4(z), F_8(z), F_{16}(z)$ we substitute the representations of $A_2,\ldots,A_7$ and $m$
in terms of $a_1,\ldots,a_7$. 
We extract the gcd-s $W_4,W_8,W_{16}$ of the coefficients of $F_4(z), F_8(z), F_{16}(z)$, 
respectively, then 
we obtain
\[
F_4(z)=W_4\cdot G_4(z),\;F_8(z)=W_8\cdot G_8(z),\;F_{16}(z)=W_{16}\cdot G_{16}(z),
\]
where $G_4(z),G_8(z),G_{16}(z)$ are polynomials of $z_1,\ldots,z_7$ with integer coefficients,
and
\[
W_4W_8W_{16}=\frac{A_7^{28}}{A_2A_3A_4A_5A_6A_7}.
\]
Therefore
\[
\ind(\vartheta)=G_4(z)\cdot G_8(z)\cdot G_{16}(z).
\]
Calculating explicitly $G_{16}(z)-a_2^2a_6^2G_8^2(z)$ we obtain 
\[
8a_1a_3a_5a_7 | (G_{16}(z)-a_2^2a_6^2G_8^2(z)).
\]
If $K$ is monogenic, then for some $z_1,\ldots,z_7\in\Z$ we have 
 $\ind(\vartheta)=1$, hence $G_4(z)=\pm 1, G_8(z)=\pm 1, G_{16}(z)=\pm 1$,
that is the above divisibility relation implies 
\[
8a_1a_3a_5a_7|(a_2^2a_6^2\pm 1).
\]
with at least one of the signs.
\hfill$\Box$

\vspace{0.5cm}

\vspace{0.5cm}
{The existence of prime common index divisors was first established in $1871$ by Dedekind who exhibited examples in cubic number fields. For example, he considered the cubic field $K$ generated  by a root of $x^3-x^2-2x-8$ and he showed that the prime $2$ splits completely in $\Z_K$. So, if we suppose that $K$ is monogenic, then we would be able to find a cubic polynomial generating $K$, that splits completely into distinct polynomials of degree $1$ in $\mathbb{F}_2[x]$. Since there are only $2$ distinct polynomials of degree $1$ in $\mathbb{F}_2[x]$, this is impossible. Based on these ideas and using Kronecker's theory of algebraic numbers, Hensel  gave a necessary and sufficient condition on the so-called "index divisors" for any prime integer  $p$ to be  a prime common index divisor \cite{He2}.}\\
	
In order to prove  Theorem \ref{npib},  we need the following two lemmas. The  first one is an immediate consequence of {Dedekind's} theorem. The second one follows from  \cite[Corollary 3.8]{GMN}. 

\vspace{0.5cm}

\begin{lem} \label{comindex}
 Let  $p$ be  a rational prime integer and $K$  a number field. For every positive integer $f$, let $P_f$ be the number of distinct prime ideals of $\Z_K$ lying above $p$ with residue degree $f$ and $N_f$  the number of monic irreducible polynomials of  $\F_p[x]$ of degree $f$.
{ If $ P_f > N_f$ for some
positive integer $f$}, then $p$ is a common index divisor  of $K$.
\end{lem}

\vspace{0.5cm}

By Hensel's correspondence and \cite[Corollary 3.3]{GMN}, the following Lemma allows to
calculate the ramification index and the residue degree of each prime ideal of $K$ lying
above a prime integer $p$.

\vspace{0.5cm}

\begin{lem} \label{residue}
 Let $p$ be a prime integer, $f(x)\in \Z_p[x]$ a monic  polynomial such that $\ol{f(x)}$ is a  power of $\ol{\ph(x)}$ for some monic polynomial $\ph\in \Z_p[x]$, whose reduction is irreducible over $\F_p$, $N_i({f})=S_i$ has a single side of slope $-\la_i$ and $R_i(f)(y)=\psi_i^{a_i}(y)$ for some monic irreducible polynomial $\psi_i\in \F_i[y]$ for every $i=1,\dots, r-1$. Let $e_i$ be the smallest positive integer satisfying $e_i\la_i\in \Z$ for every $i=1,\dots,r-1$. If $N_r(f)=T$ has a single side of degree $1$, then $f(x)$ is irreducible over $\Q_p$. Let $\p$ be the unique prime ideal of $\Q(\beta)$ lying above $p$, where $\beta$ is a root of $f(x)$. Then $e(\p/p)=e_1\times\dots\times e_r$ is  the ramification index of $\p$  and $f(\p/p)=\mbox{deg}(\ph)\times \prod_{i=1}^{r-1}\mbox{deg}(\psi_i)$ is its residue degree.
\end{lem}

\vspace{0.5cm}

\vspace{0.5cm}

\noindent
{\it Proof of Theorem \ref{npib}}
\begin{enumerate}
\item  
 If $m\equiv 1 \md{32}$, then  $\ol{f(x)}=\ol{\ph(x)}^8\md2$ and $\nph{f}=S_{1}+S_{2}+S_3+S_4$ has $4$ sides  joining $(0,v)$, $(1,3)$, $(2,2)$, $(4,1)$, and $(8,0)$, where $v=\nu_2(1-m)$. Thus each side is of  degree $1$, and so  there are $4$ prime ideals of $\Z_K$ lying above $2$. Since there are only $2$ monic irreducible polynomial of degree modulo $2$. By Lemma \ref{comindex}, $2$ divides $i(K)$, and so $K$ is not monogenic.
  \item  
 If $\nu_2(16-m)\ge 8$; $m\equiv 272\md{512}$, then $N_1(f)=S$ has a single side of slope $-\la_1=-1/2$ (see Figure 2). Then $e_1=1$ is the ramification index of $S$. Also $N_2(f)=T_1+T_2+T_3$  has $3$ sides joining the points $(0,v)$, $(1,12)$, $(2,10)$, and $(2,8)$ (see Figue 3) with 
  $v=\om_2(16-m)\ge 16$, residual polynomials $R_{i2}(f)(y)=y+1$ for $i=1,2$, and $R_{32}(f)(y)=(y+1)^2$. By Theorem \ref{highpol}, $f(x)=f_1(x)\times f_2(x)\times f_3(x)$ in $\Z_2[x]$ where $f_i(x)$ is monic irreducible over $\Q_2$ for $i=1,2$ and $f_3(x)$ is monic. In order to complete the factorization of $f_3(x)$, we have to use third order Newton polygon. For $\ph_3=\ph_2+4x=x^2+4x+2$, we get  $N_3^-(f)=T$ has a single side of degree  $1$ (see Figure 4). Thus $f_3(x)$ is  irreducible over $\Q_2$. By Hensel's correspondence, there are exactly $3$ prime ideals $\p_1,\p_2,\p_3$ of $\Z_K$ lying above $2$ (each $\p_i$ is associated to the factor $f_i$). By Lemma \ref{residue},
  $e(\p_i)=e_1\times e_2=2$ for $i=1,2$ and $e(\p_3)=e_1\times e_2\times e_3=4$. As $\sum_{i=1}^3e(\p_i)f(\p_i)=8$, we conclude that the residue degree $f(\p_i)=1$ for every $i=1,2,3$.  As there are only $2$ monic irreducible polynomial in $\F_0[x]$, $2$ divides $i(K)$, and so $K$ is not monogenic.
 \begin{figure}[htbp] 
\centering

\begin{tikzpicture}[x=1cm,y=0.5cm]
\draw[latex-latex] (0,6) -- (0,0) -- (10,0) ;

\draw[thick] (0,0) -- (-0.5,0);
\draw[thick] (0,0) -- (0,-0.5);

\node at (0,0) [below left,blue]{\footnotesize $0$};

\draw[thick] plot coordinates{(0,4)  (8,0)};

\node at (4,3) [above  ,blue]{\footnotesize $S$};
\end{tikzpicture}
\caption{    \large  $N_1({f})$.}
\end{figure}

\begin{figure}[htbp] 
\centering

\begin{tikzpicture}[x=1cm,y=0.5cm]
\draw[latex-latex] (0,6) -- (0,0) -- (10,0) ;

\draw[thick] (0,0) -- (-0.5,0);
\draw[thick] (0,0) -- (0,-0.5);

\node at (0,1.5) [below left,blue]{\footnotesize $8$};
\node at (4,0) [below ,blue]{\footnotesize $4$};
\node at (2,0) [below ,blue]{\footnotesize $2$};
\node at (1,0) [below ,blue]{\footnotesize $1$};
\draw[thick] plot coordinates{(0,1) (6,1)};
\draw[thick] plot coordinates{(0,5.5) (1,3) (2,2) (4,1)};

\node at (0.7,5) [above  ,blue]{\footnotesize $T_{1}$};
\node at (1.5,2.7) [above   ,blue]{\footnotesize $T_{2}$};
\node at (3,1.8) [above   ,blue]{\footnotesize $T_{3}$};
\end{tikzpicture}
\caption{    \large  $N_2({f})$.}
\end{figure}

\begin{figure}[htbp] 
\centering

\begin{tikzpicture}[x=1cm,y=0.5cm]
\draw[latex-latex] (0,6) -- (0,0) -- (10,0) ;

\draw[thick] (0,0) -- (-0.5,0);
\draw[thick] (0,0) -- (0,-0.5);

\node at (0,2.5) [below left,blue]{\footnotesize $12$};
\node at (0,3.5) [below left,blue]{\footnotesize $13$};\node at (4,0) [below ,blue]{\footnotesize $4$};
\node at (2,0) [below ,blue]{\footnotesize $2$};
\node at (1,0) [below ,blue]{\footnotesize $1$};
\draw[thick] plot coordinates{(0,2) (6,2)};
\draw[thick] plot coordinates{(0,3)  (2,2)};

\node at (0.7,5) [above  ,blue]{\footnotesize $T$};
\end{tikzpicture}
\caption{    \large  $N_3^-({f})$.}
\end{figure}
\item
 If $\nu_2(m)$ is odd, then by Theorem \ref{index},
\[
8a_1a_3a_5a_7|(a_2^2a_6^2\pm 1)
\]
 is a necessary condition for monogenity of $K$ with at least one of the signs.  
Thus if $a_2a_6\equiv 2\md8$ or $a_2a_6\equiv 6\md8$, then $8a_1a_3a_5a_7$  does not divide $(a_2a_6)^2\pm 1$
with neither of the signs. Thus $K$ is not monogenic.
\end{enumerate}
\hfill$\Box$

\noindent
{\it Proof of of Theorem \ref{mono}}\\
Assume that  $m=a^t$ with $a\neq \pm 1$  is a square free rational integer and $t\in \{3,5,7\}$. Since gcd$(t,8)=1$, let $(u,v)$ be two non-negative integers solution of  $tu-8v=1$ and $u<8$ and  $\theta=\frac{\al^u}{a^v}$. Since $\theta^8=\frac{\al^{8u}}{a^{8v}}=\frac{m^u}{a^{8v}}=\frac{a^{tu}}{a^{8v}}=a$ and $x^8-a$ is irreducible over $\Q$ (because it is an Eisenstein polynomial), hence $\Q(\theta)=\Q(\al)=K$. Thus $K$ is the octic number field generated by a root 
$\theta$ of  $g(x)=x^8-a$ with a square free integer $a$. Applying \cite[Theorem 8]{GR17}, we get the following results:
\begin{enumerate}
\item 
If  $a\not\equiv
1 \md{4}$, then $K$ is monogenic
and $\Z_K$ is generated by  $\theta$.
\item
If  $a\equiv
1 \md{4}$, then $K$ is not monogenic, except for $a=-3$.
\end{enumerate}
\hfill$\Box$

\section*{Acknowledgments}
The authors are deeply  grateful to the anonymous referee whose valuable comments and suggestions have tremendously improved the quality of this paper. Also, the  first author is deeply grateful to Professor Enric Nart who introduced him to Newton polygon techniques when he was a post-doc at CRM of Barcelona (2007--2008).


\begin{thebibliography}{99}

\bibitem{AN} {\sc S. Ahmad, T. Nakahara, and A. Hameed}, {\it On certain pure sextic fields related to a problem of Hasse}, Int. J. Alg. and Comput. 26(3) (2016), 577--583.

\bibitem{NH} {\sc S. Ahmad, T. Nakahara and S. M. Husnine}, {\it Non-monogenesis of a family
of pure sextic fields}, Arch. Sci. (Geneva) 65(7) (2012), 42.

 \bibitem{AN6} {\sc S. Ahmad, T. Nakahara, and S. M. Husnine}, {\it Power integral bases for certain pure sextic
fields}, Int. J. of Number Theory 10(8) (2014), 2257-- 2265.

\bibitem{maple}
{\sc B. W. Char, K. O. Geddes, G. H. Gonnet, M. B. Monagan, S. M. Watt} (eds.),
{\it MAPLE, Reference Manual}, Watcom Publications, Waterloo, Canada, 1988.


\bibitem{Co} {\sc H. Cohen}, {\it A Course in Computational Algebraic Number Theory}, GTM 138, Springer-Verlag Berlin  Heidelberg (1993).

\bibitem{R} {\sc R. Dedekind},  {\it \"Uber den Zusammenhang zwischen der Theorie der Ideale und der Theorie der h\"oheren Kongruenzen}, G\"ottingen Abhandlungen {\bf 23} (1878), 1--23.

\bibitem{El} {\sc L. El Fadil}, {\it On Newton polygon techniques and factorization of polynomial over Henselian valued fields},  J. of Algebra and  its Appl. {\bf 19}(10) (2020), 2050188.

\bibitem {E6} {\sc L. El Fadil}, {\it On Power integral bases for certain pure sextic fields},
Bol. Soc. Paran. Math. to appear,  DOI: 10.5269/bspm.42373 .

\bibitem{EMN}{\sc L. El Fadil, J. Montes E. Nart}, {\it Newton polygons and $p$-integral bases of quartic number fields}, J. Algebra and Appl. {\bf 11}(4)(2012), Article ID 1250073.

\bibitem{gaalbook} {\sc I.Ga\'al}, {\it Diophantine equations and power integral bases, Theory and algorithms}, 
Second edition, Boston, Birkh\"auser, 2019.

\bibitem{GR14} {\sc  I. Ga\'al and L. Remete}, 
{\it Binomial Thue equations and power integral bases in pure quartic fields},
JP J. Algebra Number Theory Appl. {\bf 32} (2014), No. 1, 49--61.

\bibitem{GR17} {\sc  I. Ga\'al and L. Remete}, {\it Power integral bases and monogenity of pure  fields}, J. of Number Theory {\bf 173} (2017) 129--146.
\bibitem{GMN}{\sc J. Guardia, J. Montes and E. Nart}, {\it Newton
polygons of higher order in algebraic number theory}, Trans. Amer. Math. Soc. {\bf 364} (1) (2012) 361--416.
\bibitem{GN}
{ {\sc J. Guardia and  E. Nart}, {\it Genetics of polynomials over local fields}, 
Contemp. Math. {\bf 637} (2015), 207--241.   }

\bibitem{HN15}
{\sc A. Hameed and T. Nakahara}, {\it Integral bases and relative monogeneity of pure octic fields},
Bull. Math. Soc. Sci. Math. R\'epub. Soc. Roum., {\bf 58(106)} (2015),  No. 4, 419-433.

\bibitem{hasse63} {\sc H. Hasse}, {\it Zahlentheorie}, Akademie-Verlag, Berlin, 1963.

\bibitem{He2} {\sc K. Hensel}, 
{\it Arithemetishe untersuchungen uber die gemeinsamen ausserwesentliche Discriminantentheiler einer Gattung}, 
J. Reine Angew  Math, {\bf 113} (1894), 128--160.

\bibitem{hensel08} {\sc K. Hensel}, {\it Theorie der algebraischen Zahlen}, Teubner Verlag, Leipzig, Berlin, 1908.

\bibitem{MN}  {\sc J. Montes and E. Nart}, {\it On a theorem of Ore}, J.  Algebra {\bf 146}(2) (1992), 318--334.

\bibitem{Mc} {\sc S. MacLane}, {\it A construction for absolute values in polynomial rings},
 Trans. Amer. Math. Soc. 	{\bf 40} (1936), 363--395.

\bibitem{Neu} {\sc J.  Neukirch}, {\it Algebraic Number Theory},
Springer-Verlag, Berlin, 1999.

\bibitem{O}{\sc O. Ore}, {\it Newtonsche Polygone in der Theorie der algebraischen Korper}, 
Math. Ann., {\bf 99} (1928), 84--117.

\end{thebibliography}
  \end{document}